\documentclass[12pt]{article}

\usepackage{authblk}
\usepackage{ulem}
\usepackage{amsmath,amsfonts,amssymb,lmodern,geometry,enumerate}
\usepackage[font=small,labelfont=bf]{caption}
\usepackage{tkz-euclide}
\usepackage{subcaption}
\usepackage[T1]{fontenc}
\usepackage[english]{babel}
\usepackage{lmodern}
\usepackage{scalefnt}
\usepackage{dsfont}
\usepackage{stmaryrd}
\usepackage{color}
\usepackage{bm,bbm}
\usepackage{mathrsfs,url,color}
\usepackage{breakcites} %for breaking the citations
\usepackage{textcase}
\usepackage{graphicx}
\usepackage{wrapfig}
\usepackage{flushend,cuted}
\usepackage{bm}
\usepackage{tabularx}
\usepackage{color}
\usepackage{indentfirst}
\usepackage{amssymb}
\usepackage{xparse}
\usepackage{tikz}
\usepackage{pgfplots}

\pgfplotsset{compat=newest}
\usepackage{mdwlist}
\usepackage{amsmath}
\usepackage{amsthm}
\usepackage{nameref}
\usetikzlibrary[intersections,
positioning,
petri,
backgrounds,
fit,
decorations.pathmorphing,
arrows,
arrows.meta,
bending,
calc,
intersections,
through,
backgrounds,
shapes.geometric,
quotes,
matrix,
trees,
shapes.symbols,
graphs,
math,
patterns,
external,
scopes,
matrix,
lindenmayersystems,
shapes.callouts,
shapes.misc,
angles,
shapes.arrows,
shadings]

\makeatletter
\def\@noindentfalse{\global\let\if@noindent\iffalse}
\def\@noindenttrue {\global\let\if@noindent\iftrue}
\def\@aftertheorem{%
  \@noindenttrue
  \everypar{%
    \if@noindent%
      \@noindentfalse\clubpenalty\@M\setbox\z@\lastbox%
    \else%
      \clubpenalty \@clubpenalty\everypar{}%
    \fi}}

\newtheorem{thm}{Theorem}[section]
\AfterEndEnvironment{thm}{\@aftertheorem}

\AfterEndEnvironment{prop}{\@aftertheorem}
\newtheorem{defi}[thm]{Definition}
\AfterEndEnvironment{defi}{\@aftertheorem}
\newtheorem{lma}[thm]{Lemma}
\AfterEndEnvironment{lma}{\@aftertheorem}
\newtheorem{cor}[thm]{Corollary}
\AfterEndEnvironment{cor}{\@aftertheorem}
\newtheorem{re}[thm]{Remark}
\AfterEndEnvironment{re}{\@aftertheorem}
\newtheorem{ap}{Application}

\theoremstyle{remark}
\newcommand{\diam}{{\rm diam}}
\newcommand{\law}{\mathscr{L}}

\newcommand{\Pro}{\mathbb{P}} 
\newcommand{\prob}{\Pro}
\newcommand{\E}{\mathbb{E}}
\newcommand{\mean}{\E}

\newcommand{\Var}{\mathrm{Var}}
\newcommand{\Vol}{\mathrm{Vol}}

\def\d{{\delta}}
\def\e{{\epsilon}}
\def\iid{{i.i.d.}}

\def\[{\left[}
\def\]{\right]}
\def\({\left(}
\def\){\right)}

\def\var{{\rm Var}}

\def\ignore#1{}
\def\Ref#1{(\ref{#1})}
\def\qed{\hfill\hbox{${\vcenter{\vbox{
					\hrule height 0.4pt\hbox{\vrule width 0.4pt height 6pt
						\kern5pt\vrule width 0.4pt}\hrule height 0.4pt}}}$}}

\newcounter{con}%for positive constants
\stepcounter{con}
\newcommand{\qcon}[1]{\addtocounter{con}{1}}
\newcounter{cproofa}%for positive constants

\newcounter{cproofb}%for positive constants

\newcounter{cproofc}%for positive constants
\newcounter{cproofd}
\newcounter{cproofe}
\newcounter{cprooff}
\newcounter{cproofg}

\newcounter{xiaa}
\newcounter{xiab}
\newcounter{xiac}
\newcounter{xiad}
\newcounter{xiae}
\newcounter{xiaf}
\newcounter{xiag}

\newcounter{xiaha}
\newcounter{xiahb}
\newcounter{xiahc}
\newcounter{xiahd}
\newcounter{xiaia}
\newcounter{xiaib}
\newcounter{xiaic}
\newcounter{xiaid}
\newcounter{xiaie}

\numberwithin{equation}{section}

\ignore{\makeatletter
	\renewcommand\section{\@startsection {section}{1}{\z@}%
		{-3.5ex \@plus -1ex \@minus -.2ex}%
		{1.3ex \@plus.2ex}%
		{\center\small\sc\MakeTextUppercase}}
	\def\subsection#1{\@startsection {subsection}{2}{0pt}%
		{-3.5ex \@plus -1ex \@minus -.2ex}%
		{1ex \@plus.2ex}%
		{\bf\mathversion{bold}}{#1}}
	\def\subsubsection#1{\@startsection{subsubsection}{3}{0pt}%
		{\medskipamount}%
		{-10pt}%
		{\normalsize\itshape}{\kern-2.2ex. #1.}}
	\makeatother}
\allowdisplaybreaks[4]

\begin{document}
	
\title{\sc\bf\large\MakeUppercase{
			Convergence rate for geometric statistics of point processes with fast decay dependence
}}

\author[1]{Tianshu Cong\thanks{{\sf{email: tcong1@student.unimelb.edu.au.}} Work supported by a Research Training Program Scholarship, a faculty of science postgraduate writing-up award and a Xing Lei Cross-Disciplinary PhD Scholarship in Mathematics and Statistics at the University of Melbourne..}}
\author[1]{Aihua Xia\thanks{{\sf{email: aihuaxia@unimelb.edu.au}}. Work supported by the Australian Research Council Grant No DP190100613.}}

%\affil[1]{%
\affil[1]{%
School of Mathematics and Statistics, the University of Melbourne, Parkville VIC 3010, Australia} 

\date{\today}
	
\maketitle
\vskip-1cm
\begin{abstract}
		\cite{BYY19} established central limit theorems for geometric statistics of point processes having fast decay dependence. As limit theorems are of limited use unless we understand their errors involved in the approximation, in this paper, we consider the rates of a normal approximation in terms of the Wasserstein distance for statistics of point processes on $\mathbb{R}^d$ satisfying fast decay dependence. We demonstrate the use of the theorems for statistics arising from two families of point processes: the rarified Gibbs point processes and the determinantal point processes with fast decay kernels.
\end{abstract}

	
	
\section{Introduction}
	
Random events in space and time can be represented as point processes and statistics arising from such random events are often of the form $\sum_{x\in \Xi\cap A} \eta(x,\Xi)$, where $\Xi$ is a point process on $\mathbb{R}^d$, $A\subset \mathbb{R}^d$ is a bounded Borel set, and $\eta(x,\Xi)$, known as a {\it score function}, is a functional representing the interaction between the point at $x$ and the point process $\Xi$. The form dates back to \cite{BHH59,Steele81} with $\Xi$ as a binomial point process and \cite{AB93} with $\Xi$ being a Poisson point process. As the density of points or the size of the observation window $A$  increases, various limit theorems of the functionals have been established since then, see, e.g., \cite{S12,S16,LSY19} and references therein. If $\eta(x,\Xi)=\bm{1}_{x\in \Xi}$ and $\Xi$ is an $\alpha-$determinantal point process with some condition on $\alpha$, \cite{S02,ST03,NS12} proved that the counting statistics are asymptotically normal. The asymptotic normality also holds when the input process $\Xi$ is a Gibbsian point process \cite{SY13,XY15}.  The common feature leading to  these limit theorems is the local dependence \cite[p.~133]{CS04} in the sense that each summand is affected by a small neighbourhood only and hence it makes a nearly independent contribution to the statistics of interest.  More generally, when the underline point process $\Xi$ is a stationary point process and has fast decay dependence, and the score function $\eta(x,\Xi)$ is determined by points of $\Xi$ not too far away from $x$, it is possible to establish central limit theorems for such statistics \cite{BYY19}. In this paper, we aim to quantify the errors associated with the limit theorems in \cite{BYY19} because the limit theorems are of limited practical value unless we understand the magnitude of the errors involved in the approximation of these statistics. 

The local influence of the score function $\eta$ can be controlled through
the concept of stabilisation \cite{BX,PY01,PY03,PY05,Pe07a,Pe07b,SY13,XY15} and the fast decay of dependence of the input point process can be represented through its 
correlation functions \cite[(1.2)]{LSY19} being controlled by a correlation decay function \cite[Definition~1.2]{LSY19}. The limit theorems in \cite{LSY19} are formulated in terms of the Palm distributions of all orders \cite[p.~110]{Kallenberg83} of the point process, and their proofs hinge on the Marcinkiewicz theorem \cite[Lemma~3]{S02} and\ignore{ on} the method of cumulants. The Marcinkiewicz theorem is a handy tool to prove the central limit theorems, but when we aim for the errors of approximation, it seems impractical to use the cumulants to control the errors because these quantities are hard to obtain in applications. For this reason, we impose the condition of the decay of dependence through the $\beta$-mixing coefficient in {\it A2.0 Exponentially Decay Dependence} (EDD). 
We will show in Section~\ref{PPA2.0} that the Gibbs point processes with nearly finite range potentials, a class of the determinantal point processes with fast decay kernels, the $r$-dependent point processes and the Boolean models all possess the EDD property. We also consider a wider class of score functions with different types of edge effects. In Section~\ref{Applications}, we demonstrate the use of the theorems for statistics arising from the rarified Gibbs point processes and the determinantal point processes with fast decay kernels. For ease of reading, we postpone the proofs of the main results to Section~\ref{Theproofs}.

The asymptotic normality depends on the lower bound of the variance of the statistics, and showing the order of the variance is
itself an interesting but hard topic \cite[(1.26) and Remark~(iii)]{LSY19}. In an attempt to recover the volume order of the variances, we obtained Theorem~\ref{thmvar}, which is an analogue of \cite[Theorem~1.15]{LSY19}. 
	
	\section{General results}\label{Generalresults}
	
	%In applications, the behaviour of statistics of point processes is of considerable interest. In this paper, we concentrate on the convergence rate of point processes satisfying certain mixing conditions. 
	To start with, we recall the definition of the marked point processes on $(\mathbb{R}^d,\mathscr{B}(\mathbb{R}^d))$ with marks in a measurable space $(T,\mathscr{T})$, where $\mathbb{R}^d$ is equipped with the Euclidean norm $\|\cdot\|$ and the Borel $\sigma$-algebra $\mathscr{B}(\mathbb{R}^d)$, and $\mathscr{T}$ is a $\sigma$-algebra on $T$.
Let $\bm{S}:=\mathbb{R}^d\times T$ be equipped with the product $\sigma$-field $\mathscr{S}:=\mathscr{B}(\mathbb{R}^d)\times \mathscr{T}$. We use~$\bm{C}_{\bm{S}}$ to denote the space of all locally finite (with respect to the first coordinate in $\mathbb{R}^d$) non-negative integer-valued measures $\xi$, often called {\it configurations}, on~$\bm{S}$ such that $\xi(\{{x}\}\times T)\le 1$ for all ${x}\in\mathbb{R}^d$. The space~$\bm{C}_{\bm{S}}$ is endowed with the $\sigma$-field $\mathscr{C}_{\bm{S}}$ generated by the vague topology \cite[p.~169]{Kallenberg83}. A {\it marked point process\/} $\Xi$ on $\mathbb{R}^d$ is a measurable mapping from $(\Omega,\mathscr{F},\prob)$ to $(\bm{C}_{\bm{S}},\mathscr{C}_{\bm{S}})$ \cite[p.~49]{Kallenberg17}. The induced simple point process $\overline{\Xi}(\cdot):=\Xi(\cdot\times T)$ is called the {\it ground process\/} \cite[p.~3]{Daley08} or projection  \cite[p.~17]{Kallenberg17} of the marked point process~$\Xi$ on $\mathbb{R}^d$. We use $M_x$ to denote the mark of $\Xi$ at $x$ for $x\in \overline{\Xi}$.
	
	For a marked point process $\Xi$, let $\Xi_A$ be the restriction of $\Xi$ to $A\times T$ defined as $\Xi_A(B\times D):=\Xi((A\cap B)\times D)$, $\Xi^x$ be the shifted point process of $\Xi$ by $x$ defined as $\Xi^x(B\times D):=\Xi((B+x)\times D)$ for all $x\in \mathbb{R}^d$, $D\in \mathscr{T}$ and $A,~B\in \mathscr{B}(\mathbb{R}^d)$. We say that the marked point process $\Xi$ is {\it stationary} if $\Xi\overset{d}{=}\Xi^x$ for all $x\in \mathbb{R}^d$, where $\overset{d}{=}$ stands for `equal in distribution'. To avoid using the Palm distributions of {\it all} orders and {\it all} cumulants in the approximation bounds, in this paper, the fast decay dependence of the marked point process $\Xi$ is quantified through its {\it $\beta$-mixing coefficient} \cite{VR59,Rio17}: for 
$A_1,A_2\in \mathscr{B}(\mathbb{R}^d)$,
	$$\beta_{A_1,A_2}:=\frac12\int_{\zeta_1\in \bm{C}_{A_1\times T}}\int_{\zeta_2\in \bm{C}_{A_2\times T}}\left|\mathbb{P}\left(\Xi_{A_1}\in d\zeta_1, \Xi_{A_2}\in d\zeta_2\right)-\mathbb{P}\left(\Xi_{A_1}\in d\zeta_1\right)\mathbb{P}\left(\Xi_{A_2}\in d\zeta_2\right)\right|,$$ 
where and in the following, $\bm{C}_{A_i\times T}$ and $\mathscr{C}_{A_i\times T}$ are defined in the same way as $\bm{C}_{\bm{S}}$ and $\mathscr{C}_{\bm{S}}$ with $\bm{S}$ replaced by $A_i\times T$, $i=1,2.$

	To define the decay of dependence, we set $\diam(A):=\sup\{\|x-y\|;~x,y\in A\}$, $d(A,B):=\inf\{\|x-y\|; x\in A,~y\in B\}$ for  $A,B\in \mathscr{B}(\mathbb{R}^d)$, where we used the convention $\sup\{\emptyset\}=0$ and $\inf\{\emptyset\}=\infty$. We use $\vee$ to stand for the maximum.
	
\noindent{\it A2.0 Exponentially Decay Dependence} We say that the marked point process $\Xi$ has the {\it exponentially decay dependence} (EDD) if there exist constants $\theta_0\in \mathbb{R}_0:=[0,\infty)$, $\theta_i\in \mathbb{R}_+:=(0,\infty)$, $1\le i\le 4$, such that for any $A$, $B\in \mathscr{B}(\mathbb{R}^d)$ with $d(A,B)\ge \theta_3\ln(\diam(A)\vee \diam(B)\vee \theta_4)$,  
	\begin{equation}\label{mixing}
    			\beta_{A,B}\le \theta_1(\diam(A)^{\theta_0}\vee 1)(\diam(B)^{\theta_0}\vee 1)e^{-\theta_2 d(A,B)}.
    		\end{equation} 
		
	The idea of the EDD is that the total variation distance between the law of $(\Xi_{A},\Xi_{B})$ and the law of the independent union of $\Xi_{A}$ and $\Xi_{B}$ decays exponentially fast as the distance between $A$ and $B$ becomes large. 
    
    The following lemma says that, in applications, it is sometimes more convenient to verify the EDD via the volumes of the sets. To this end, let $\Vol (A)$ denote the volume of the set $A\in \mathscr{B}(\mathbb{R}^d)$.
    
    \begin{lma}\label{mixing.re} If there exist constants $\theta_0'\in \mathbb{R}_0$, $\theta_i'\in \mathbb{R}_+$, $1\le i\le 4$, such that for any $A$, $B\in \mathscr{B}(\mathbb{R}^d)$ with $d(A,B)\ge \theta_3'\ln(\Vol(A)\vee \Vol(B)\vee \theta_4')$, 
		$$\beta_{A,B}\le \theta_1'(\Vol(A)^{\theta_0'}\vee 1)(\Vol (B)^{\theta_0'}\vee 1) e^{-\theta_2' d(A,B)},
	$$    then $\Xi$ satisfies the EDD.
   \end{lma}
    
    \noindent{\it Proof.} Given the dimension $d$, $\text{Vol}(A)\le \frac{\pi^{d/2}}{2^d\Gamma(\frac{d}{2}+1)}\diam(A)^d$, hence 
   \Ref{mixing} follows immediately. \qed

\begin{re}\label{rebeta} The constants $\theta_4$ in the definition of the EDD and $\theta_4'$ in Lemma~\ref{mixing.re} are not essential and they can be replaced by any positive constants because the definition of $\beta$-mixing coefficient ensures that $\beta$ is non-decreasing in the sense of inclusion, i.e., $\beta_{A,B}\le \beta_{A',B'}$ for all $A,A',B,B'\in \mathscr{B}(\mathbb{R}^d)$ such that $A\subset A'$ and $B\subset B'$.
\end{re}

    Let $\Xi$ be a stationary marked point process on $\bm{S}$ with independent and identically distributed ({\it $\iid$}) marks that satisfies the EDD. Writing the law of $\Xi$ as $\mathscr{P}$ and the law of the independent marks as $\law_{T}$. The functionals we study in the paper are defined on $\Gamma_\alpha:=\[-\frac{1}{2}\alpha^{\frac{1}{d}}, \frac{1}{2}\alpha^{\frac{1}{d}}\]^d$, the cube with volume $\alpha$ on $\mathbb{R}^d$ centred at $\bm{0}$, having the forms
    $$W_\alpha:={\sum_{(x,m)\in\Xi_{\Gamma_\alpha}}\eta(\left(x,m\right), \Xi)}$$
    and
    $$\bar{W}_\alpha:={\sum_{{(x,m)\in\Xi_{\Gamma_\alpha}}}\eta(\left(x,m\right), \Xi_{\Gamma_\alpha},\Gamma_\alpha)=\sum_{{(x,m)\in\Xi_{\Gamma_\alpha}}}\eta(\left(x,m\right), \Xi,\Gamma_\alpha)}.$$
	The function $\eta$ is called a {\it score function} {(resp. {\it restricted score function})}, i.e., a measurable function on $\(\bm{S}\times \bm{C}_{\bm{S}},\mathscr{S} \times  \mathscr{C}_{\bm{S}}\)$ to $\(\mathbb{R},\mathscr{B}\(\mathbb{R}\)\)$ (resp. a function mapping $\bm{S}\times \bm{C}_{\Gamma_\alpha\times T}\times\mathscr{B}(\mathbb{R}^d)$ to $\mathbb{R}$ which is 
	 $\(\bm{S}\times \bm{C}_{\Gamma_\alpha\times T},\mathscr{S} \times  \mathscr{C}_{\Gamma_\alpha\times T}\)$ to $\(\mathbb{R},\mathscr{B}\(\mathbb{R}\)\)$  measurable when the third coordinate is fixed) and it represents the interaction between a point with its mark and the configuration of the point process. The class of score functions considered here is broader than that considered in \cite{BYY19}. More precisely, if the score function for the restricted case does not depend on the third argument, it reduces to that in \cite{BYY19}. Because the interest is in the values of the score function of the points in a configuration, for convenience, $\eta\((x,m),\mathscr{X}\){~(\mbox{resp.}~\eta\((x,m),\mathscr{X},\Gamma_\alpha\))}$ is understood as $0$ for all $x\in \mathbb{R}^d$ and $\mathscr{X}\in \bm{C}_{\bm{S}}$ such that $(x,m)\notin \mathscr{X}$. 
	
	We need Palm processes and reduced Palm processes as the tools for stating the conditions and constructing proofs. For ease of reading, we briefly recall their definitions. Let $H$ be a Polish space with Borel $\sigma$-algebra $\mathscr{B}\(H\)$ and configuration space $\(\bm{C}_H, \mathscr{C}_H\)$, let $\Upsilon$ be a point process on $\(H, \mathscr{B}\(H\)\)$ and write the mean measure of $\Upsilon$ as $\mathbb{E} \Upsilon$. The point processes $\{ \Upsilon_x : x \in H\}$ are said to be the {\it reduced Palm processes} associated with $\Upsilon$ if for any measurable function $f : (H\times \bm{C}_H,\mathscr{B}\(H\)\times \mathscr{C}_H) \rightarrow (\mathbb{R}_0,\mathscr{B}(\mathbb{R}_0))$,
	\begin{equation}\label{palm1}
		\mathbb{E} \[ \int_{H} f(x,\Upsilon)\Upsilon(dx) \] = \int_{H} \mathbb{E} f(x,\Upsilon_x+\delta_x)  \mathbb{E} \Upsilon (d x) ,
	\end{equation}
	\cite[\S~10.1]{Kallenberg83}, where $\delta_x$ is the Dirac measure at $x$. The distributions of $\Upsilon_x$ and $\Upsilon_x+\delta_x$ are respectively called the reduced Palm distribution and the Palm distribution of $\Upsilon$ at $x$. When the point process $\Upsilon$ is {\it simple}, i.e., $\mathbb{P}(\Upsilon(\{x\})\in\{0,1\}\mbox{ for all }x\in H)=1$, the Palm distribution of $\Upsilon_x+\delta_x$ can be interpreted as the conditional distribution of $\Upsilon$ given $\Upsilon(\{x\})=1$~\cite[\S~10.1]{Kallenberg83}. 
	
For the marked point process $\Xi\sim \mathscr{P}$, since the marks are independent of each other and independent of the ground process, we can adapt \Ref{palm1} to the marked point process $\Xi$. For $f: \(\bm{S}\times \bm{C}_{\bm{S}},\mathscr{S}\times \mathscr{C}_{\bm{S}}\)\to (\mathbb{R}_0,\mathscr{B}(\mathbb{R}_0))$, recalling that $M_x$ is the mark of $\Xi$ at the point $x\in\overline{\Xi}$, we have
	\begin{equation}\label{palm4}
		\mathbb{E} \[ \int_{\mathbb{R}^d} f((x,M_x),\Xi)\overline{\Xi}(dx) \] = \int_{\mathbb{R}^d}  \mathbb{E} f((x,M),\Xi_x+\d_{(x,M)})  \lambda dx,
	\end{equation}
    where $M\sim\mathscr{L}_{T}$, $\lambda:=\mathbb{E}\(\overline{\Xi}(d\bm{0})\)/d\bm{0}$, $\Xi_x$ is the point process obtained by attaching the reduced Palm process $\overline{\Xi}_x$ of $\overline{\Xi}$ with \iid\ marks following $\law_T$, and $M$ is independent of $\Xi_x$. Without loss, we use the convention that $M_x$ is independent of $\Xi_x$ throughout the paper.  Hence, when we need to emphasise the location of the mark, we can replace $M$ by $M_x$ at the right hand side of \Ref{palm4}.
	
	The following assumptions are adapted from those in \cite{CX20} which were initiated in \cite{PY01} and further refined in \cite{XY15}.
	
	\newpage
	\noindent{\it A2.1 Stabilisation}
	
	For a locally finite configuration $\mathscr{X}$ and $z\in \bm{S}\cup\{\emptyset\}$, write $\mathscr{X}^{\lbag z \rbag}=\mathscr{X}$ if $z=\emptyset$ and $\mathscr{X}^{\lbag z \rbag}=\mathscr{X}\cup\{z\}$ otherwise. We use $B(x,r)$ to stand for the ball with centre $x$ and radius $r\ge 0$. %The notion of stabilisation is introduced in \cite{PY01}, and we adapt it to our setup as follows. 
	\begin{defi}\label{defi4} (unrestricted case)
		A score function $\eta$ on $\bm{S}$ is range-bound (resp. exponentially stabilising) with respect to $\mathscr{P}$ if for all $x\in \mathbb{R}^d$, $z\in \bm{S}\cup\{\emptyset\}$, and almost all realisations $\mathscr{X}$ of the marked point process $\Xi_x$, there exists a radius of stabilisation $R:=R(x):=R((x,m_x),\mathscr{X}^{\lbag z \rbag})\in(0,\infty)$, such that for all locally finite $\mathscr{Y}\subset (\mathbb{R}^d\backslash B(x,R))\times T$, we have 
		$$
		\eta\left(\left(x,m_x\right), \left[\mathscr{X}^{\lbag z \rbag}\cap\left(B(x,R)\times T\right)\right]\cup \mathscr{Y} \right)=\eta\left(\left(x,m_x\right),\mathscr{X}^{\lbag z \rbag}\cap\left(B(x,R)\times T\right) \right)
		$$
		and the tail probability $$\tau(t):=\sup_{(x,m_x)\in \mathbb{R}^d\times \rm{supp}(\mathscr{L}_T)}\sup_{z\in \mathbf{S}\cup\{\emptyset\}} \mathbb{P}\left(R((x,m_x),\Xi_x^{\lbag z \rbag}+\delta_{(x,m_x)})\ge t\right)$$ satisfies $
		\tau(t)=0\mbox{ for some }t>0 
		$
		({resp. }$\tau(t)\le C_1e^{-C_2t}{~\mbox{for all }t>0}$, where $C_1$ and $C_2$ are positive constants independent of $t$). 
	\end{defi}
	
	For the functionals with the input of a restricted marked point process, we have the following counterpart of stabilisation. Note that the score function for the restricted input is not affected by points outside $\Gamma_\alpha$.
		
	\begin{defi}\label{defi4r} (restricted case) We say that the score function $\eta$ is range-bound (resp. exponentially stabilising) with respect to $\mathscr{P}$ if for $\alpha\in \mathbb{R}_+$, $x\in \Gamma_\alpha$, and $z\in (\Gamma_\alpha \times T)\cup\{\emptyset\}$,  almost all realisations $\mathscr{X}$ of the marked point process  $\Xi_x$, there exists a radius of stabilisation $\bar{R}:=\bar{R}(x,\alpha):=\bar{R}((x,m_x),\alpha,\mathscr{X}^{\lbag z \rbag})\in(0,\infty)$ such that for all locally finite $\mathscr{Y}\subset (\Gamma_\alpha\backslash B(x,R))\times T$, we have 
		\begin{align*}
			&\eta\left(\left(x,m_x\right), \left[\mathscr{X}_{\Gamma_\alpha}^{{\lbag z \rbag}}\cap\left(B(x,\bar{R})\times T\right)\right]\cup \mathscr{Y},\Gamma_\alpha \right)%\nonumber\\
			=\eta\left(\left(x,m_x\right),\mathscr{X}_{\Gamma_\alpha}^{{\lbag z \rbag}}\cap\left(B(x,\bar{R})\times T\right) ,\Gamma_\alpha\right)%\label{defi4.1}
		\end{align*}
		and the tail probability $$\bar{\tau}(t):=\sup_{(x,m_x)\in \mathbb{R}^d\times \rm{supp}(\mathscr{L}_T)}\sup_{\alpha\in \mathbb{R}_+}\sup_{z\in (\Gamma_\alpha\times T)\cup\{\emptyset\}} \mathbb{P}\left(\bar{R}((x,m_x),\alpha, \Xi_x^{{\lbag z \rbag}}+\delta_{(x,m_x)})\ge t \right)$$ satisfies $\bar{\tau}(t)=0$ for some $t>0$ (resp. $\bar{\tau}(t)\le C_1e^{-C_2t}$ for all $t>0$, where $C_1$ and $C_2$
		are some positive constants independent of $t$). 
	\end{defi}

	\noindent{\it A2.2 Translation Invariance}
	
	We write $d(x,A):=\inf\{d(x,y);~y\in A\}$, $A\pm B:=\{x\pm y; \ x\in A,\ y\in B\}$ for $x\in \mathbb{R}^d$ and  $A,B\in\mathscr{B}\(\mathbb{R}^d\)$. Recall that the shift operator is defined as $\Xi^x(\cdot\times D):=\Xi((\cdot+x)\times D)$ for all $x\in \mathbb{R}^d$, $D\in \mathscr{T}$. 
	
	\newpage
	{\it A2.2.1 Unrestricted Case:} 
	
	\begin{defi}\label{invar} The score function $\eta$ is {\it translation invariant} if for all locally finite configuration $\mathscr{X}$, $x\in \mathbb{R}^d$ and $m\in T$,
		{$\eta((\bm{0},m),\mathscr{X})=\eta((x,m),\mathscr{X}^{-x})=:g(\mathscr{X})$.} 
	\end{defi}
	
	{\it A2.2.2 Restricted Case:} 
	
	As a translation may send a configuration to the outside of $\Gamma_\alpha$, resulting in a completely different configuration inside $\Gamma_\alpha$, it is necessary to focus on the part that affects the score function. Therefore, 
	we expect the score function to take the same value for two configurations if the parts within their stabilising radii are completely inside $\Gamma_\alpha$ and one is a translation of the other. More precisely, we have the following definition.
	
	\begin{defi}\label{traninvres0} A stabilising score function $\eta$ with stabilisation radius $\bar{R}$ is called {\it translation invariant} if for any $\alpha>0$, $x\in\Gamma_\alpha$  
		and $\mathscr{X}\in \bm{C}_{\bm{S}}$ such that $\bar{R}((x,m),\alpha,\mathscr{X})\le d(x,\partial\Gamma_\alpha)$, where $\partial A$ stands for the boundary of $A$, then $\eta\((x,m),\mathscr{X},\Gamma_\alpha\)=\eta\((x',m),\mathscr{X}',\Gamma_{\alpha'}\)$ and $\bar{R}((x',m),\alpha',\mathscr{X}')=\bar{R}((x,m),\alpha,\mathscr{X})$ for all $\alpha'>0$, $x'\in\Gamma_{\alpha'}$ and $\mathscr{X}'\in \bm{C}_{\bm{S}}$ such that $\bar{R}((x',m),\alpha',\mathscr{X}')\le d(x',\partial\Gamma_{\alpha'})$ and $\(\mathscr{X}'_{B(x',\bar{R}((x,m),\alpha,\mathscr{X}) )}\)^{x'}=\(\mathscr{X}_{B(x,\bar{R}((x,m),\alpha,\mathscr{X}) )}\)^{x}$.
	\end{defi}
	
	Noting that there is a tacit assumption of consistency in Definition \ref{traninvres0}, which implies that if $\eta$ is translation invariant in Definition \ref{traninvres0}, then there exists a $\bar{g}: \bm{C}_{\bm{S}}\rightarrow \mathbb{R}$ such that $$\lim_{\alpha\rightarrow \infty}\eta\((0,m),\mathscr{X},\Gamma_\alpha\)=\bar{g}\(\mathscr{X}\)$$ for $\law_{T}$ almost all $m\in T$ and almost all realisations $\mathscr{X}$ of the marked point process $\Xi\sim\mathscr{P}$. The limit $\bar g$ ensures that for each score function $\eta$ satisfying the translation-invariance in Definition \ref{traninvres0}, there exists a score function for the unrestricted case by setting $\bar{\eta}((x,m),\mathscr{X}):=\bar{g}(\mathscr{X}^{x})\mathbf{1}_{(x,m)\in \mathscr{X}}$ with the radius of stabilisation $R((x,m),\mathscr{X})=\lim_{\alpha\to\infty}\bar{R}((x,m),\alpha,\mathscr{X})$. Consequently, $\bar{\eta}$ is range bound (resp. exponentially stabilising) in the sense of Definition~\ref{defi4} if $\eta$ is range bound (resp. exponentially stabilising) in the sense of Definition~\ref{defi4r}. Moreover, if $B(x, R(x)) \subset \Gamma_\alpha$, then $\bar{R}(x, \alpha)= R(x)$, and if $B(x, R(x)) \not\subset \Gamma_\alpha$, then $\bar{R}(x,\alpha)>d(x,\partial \Gamma_\alpha)$, but there is no definite relationship between $\bar{R}$ and $R$.

	\noindent{\it A2.3 Moment condition} 
	
	We need moment conditions of both the marked point process $\Xi$ and the score function $\eta$. We say that the marked point process $\Xi\sim\mathscr{P}$ satisfies the $k$th moment condition if there exists a nonempty open set $B\subset \mathbb{R}^d$ such that \begin{equation}\label{moment0}
		\mathbb{E}\(\overline{\Xi}(B)^k\)< \infty.
	\end{equation} 
	
	For the score function $\eta$, there are two cases to consider.
	
	{\it Unrestricted Case:} The score function $\eta$ is said to satisfy the $k$th moment condition if
	\begin{equation}\label{thm2.1}
		\mathbb{E}\left[\left|\eta\left(({\bf 0},M_{\bf 0}),\Xi_{\bm{0}}+\delta_{(
			\bm{0},M_{\bf 0})}\right)\right|^k \right]<\infty.
	\end{equation}	
	
	{\it Restricted Case:} The score function $\eta$ is said to satisfy the $k$th moment condition if there exists a positive constant $C$ such that
	\begin{equation}\label{thm2.1r}
	\sup_{\alpha\in\mathbb{R}_+}\sup_{x\in\Gamma_\alpha}\mathbb{E}\left[\left|\eta\left((x,M_x),(\Xi_{x})_{\Gamma_\alpha}+\delta_{(
		x,M_x)}\right)\right|^k \right]\le C.
	\end{equation}
	
	One can verify that if $\eta$ is exponentially stabilising and satisfies \Ref{thm2.1r}, then its induced $\bar{\eta}$ satisfies the moment condition of the same order in the sense of \Ref{thm2.1}.
	
	\noindent{\it A2.4 Variation Condition} 
	
	The speed of convergence of the normal approximation is determined by the order of $\var(W_\alpha)$ or $\var(\bar{W}_\alpha)$. If we have the following condition, we can prove that the variances $\var(W_\alpha)$ and $\var(\bar{W}_\alpha)$ have the same order as the volume $\alpha$, cf. \cite[(1.22)]{BYY19}.
	
	{\it Unrestricted Case:} The score function is said to satisfy the variation condition if 
	\begin{align}
		\sigma^2=&\mathbb{E}\(\eta((\bm{0},M_{\bm{0}}), \Xi_{\bm{0}}+\delta_{(0,M_{\bm{0}})})^2\)\lambda\nonumber\\
		&+\frac{\mathbb{E}\(\int_{\mathbb{R}^d\backslash\{\bm{0}\}}\(\eta((x,M_x), \Xi)\overline{\Xi}(dx)-Pdx\)\(\eta((\bm{0},M_{\bm{0}}), \Xi)\overline{\Xi}(d\bm{0})-Pd\bm{0}\)\)}{d\bm{0}}>0,\label{non-sin}\end{align}  
	where $P:=\lambda\mathbb{E}(\eta((\bm{0},M_{\bm{0}}), \Xi_{\bm{0}}+\delta_{(0,M_{\bm{0}}})))$.
	
	{\it Restricted Case:} We define the variation condition when the score function is exponentially stabilising. The score function $\eta$ for restricted input satisfies the variation condition if it is exponentially stabilising, and the corresponding $\bar{\eta}$ satisfies
		\begin{align}
			\bar{\sigma}^2=&\mathbb{E}\(\bar{\eta}((\bm{0},M_{\bm{0}}), \Xi_{\bm{0}}+\delta_{(0,M_{\bm{0}})})^2\)\lambda\nonumber\\
			&+\frac{\mathbb{E}\(\int_{\mathbb{R}^d\backslash\{\bm{0}\}}\(\bar{\eta}((x,M_x), \Xi)\overline{\Xi}(dx)-\bar{P}dx\)\(\bar{\eta}((\bm{0},M_{\bm{0}}), \Xi)\overline{\Xi}(d\bm{0})-\bar{P}d\bm{0}\)\)}{d\bm{0}}>0,
			\label{non-sinr}
	\end{align} 
	where $\bar{P}:=\lambda\mathbb{E}(\bar{\eta}((\bm{0},M_{\bm{0}}), \Xi_{\bm{0}}+\delta_{(0,M_{\bm{0}})}))$.
	
	The above conditions are generally difficult to verify, we refer the interested readers to the discussion at Remark~(iii) of \cite[Theorem~1.14]{BYY19}. For this reason, we formulate the bounds of approximation errors in terms of the following variation conditions. If $f_1$ and $f_2$ are two functions satisfying $\liminf_{x \to \infty} f_1(x)/f_2(x) > 0$, then
we write $f_1(x) = \Omega(f_2(x))$ as $x\to\infty$.

		{\it Unrestricted Case:} $\var(W_\alpha)= \Omega(\alpha^\nu)$ for some $\nu \in (\frac{2}{3}, 1]$ as $\alpha\to\infty$.
		 
	{\it Restricted Case:} $\var(\bar{W}_\alpha)=\Omega(\alpha^\nu)$ for some $\nu \in (\frac{2}{3}, 1]$ as $\alpha\to\infty$. 
	
	In the proof of Theorem~\ref{thmvar} and Remark~\ref{re1}, we can see that under the conditions of Theorem~\ref{thmvar}, the variation condition~\Ref{non-sin} (resp. \Ref{non-sinr}) holds if and only if the above condition holds for $\nu=1$.
		
Given the variation conditions, we can establish the convergence rate in terms of the Wasserstein distance defined as \begin{equation}d_W(X,Y):=\sup_{h\in \mathscr{F}_{\rm Lip}}\mathbb{E}\(h(X)-h(Y)\),\label{defWass}
%=\inf_{\tilde{X}\overset{d}{=} X,\tilde{Y}\overset{d}{=} Y}\mathbb{E}(|\tilde{X}-\tilde{Y}|)
\end{equation}
where $\mathscr{F}_{\rm Lip}$ is the set of all Lipschitz functions $h$ on $\mathbb{R}$ such that $|h(x)-h(y)|\le |x-y|$ for all $x,y\in\mathbb{R}$. Our main result for $W_\alpha$ (unrestricted case) is summarised below.
	
    \begin{thm}\label{thm2.a1} Assume that the score function $\eta$ is translation invariant in Definition~\ref{invar} and satisfies the sixth moment condition~\Ref{thm2.1}, $\Xi$ satisfies the EDD, the fifth moment condition~\Ref{moment0} and $\var(W_\alpha)= \Omega(\alpha^\nu)$ for some $\nu \in (\frac{2}{3}, 1]$ as $\alpha\to\infty$.
    	\begin{description} 
    		\item{(i)} If $\eta$ is range-bound as in Definition~\ref{defi4}, then
    		$$d_W\(\frac{W_\alpha-\mean W_{\alpha}}{\sqrt{\var(W_{\alpha})}},Z\)\le O\(\alpha^{-\frac{3}{2}\nu+1}\).$$ 
    		\item{(ii)} If $\eta$ is exponentially stabilising as in Definition~\ref{defi4}, then
    		$$d_W\(\frac{W_\alpha-\mean W_{\alpha}}{\sqrt{\var(W_{\alpha})}},Z\)\le O\(\alpha^{-\frac{3}{2}\nu+1}\ln(\alpha)^{5d}\).$$ 
    	\end{description}
    \end{thm}

	\begin{cor}\label{thm2} Assume that the score function $\eta$ is translation invariant in Definition~\ref{invar} and satisfies the sixth moment condition~\Ref{thm2.1}, $\Xi$ satisfies the EDD, the fifth moment condition~\Ref{moment0} and \Ref{non-sin} holds.
		\begin{description} 
			\item{(i)} If $\eta$ is range-bound as in Definition~\ref{defi4}, then
			$$d_W\(\frac{W_\alpha-\mean W_{\alpha}}{\sqrt{\var(W_{\alpha})}},Z\)\le O\(\alpha^{-\frac{1}{2}}\).$$ 
			\item{(ii)} If $\eta$ is exponentially stabilising as in Definition~\ref{defi4}, then
			$$d_W\(\frac{W_\alpha-\mean W_{\alpha}}{\sqrt{\var(W_{\alpha})}},Z\)\le O\(\alpha^{-\frac{1}{2}}\ln(\alpha)^{5d}\).$$ 
		\end{description}
\end{cor}

	The restricted case is of interest in many applications, by adapting the conditions accordingly, we can show that the main result for $\bar{W}_\alpha$ (restricted case) also holds.

	 \begin{thm}\label{thm2.a2} Assume that the score function $\eta$ is translation invariant in Definition~\ref{traninvres0} and satisfies the sixth moment condition~\Ref{thm2.1r}, $\Xi$ satisfies the EDD, the fifth moment condition~\Ref{moment0} and $\var(\bar{W}_\alpha)= \Omega(\alpha^\nu)$ for some $\nu \in (\frac{2}{3}, 1]$ as $\alpha\to\infty$.
		\begin{description} 
			\item{(i)} If $\eta$ is range-bound as in Definition~\ref{defi4r}, then
			$$d_W\(\frac{\bar{W}_\alpha-\mean \bar{W}_{\alpha}}{\sqrt{\var(\bar{W}_{\alpha})}},Z\)\le O\(\alpha^{-\frac{3}{2}\nu+1}\).$$ 
			\item{(ii)} If $\eta$ is exponentially stabilising as in Definition~\ref{defi4r}, then
			$$d_W\(\frac{\bar{W}_\alpha-\mean \bar{W}_{\alpha}}{\sqrt{\var(\bar{W}_{\alpha})}},Z\)\le O\(\alpha^{-\frac{3}{2}\nu+1}\ln(\alpha)^{5d}\).$$ 
		\end{description}
	\end{thm}

 \begin{cor}\label{thm2a} Assume that the score function $\eta$ is translation invariant in Definition~\ref{traninvres0} and satisfies the sixth moment condition~\Ref{thm2.1r}, $\Xi$ satisfies the EDD, the fifth moment condition~\Ref{moment0} and 
\Ref{non-sinr} holds.
\begin{description}
    			\item{(i)} If $\eta$ is range-bound as in Definition~\ref{defi4r}, then
    			$$d_W\(\frac{\bar{W}_\alpha-\mean \bar{W}_{\alpha}}{\sqrt{\var(\bar{W}_{\alpha})}},Z\)\le O\(\alpha^{-\frac{1}{2}}\).$$
    			\item{(ii)}	If $\eta$ is exponentially stabilising as in Definition~\ref{defi4r}, then
    			$$d_W\(\frac{\bar{W}_\alpha-\mean \bar{W}_{\alpha}}{\sqrt{\var(\bar{W}_{\alpha})}},Z\)\le O\(\alpha^{-\frac{1}{2}}\ln(\alpha)^{5d}\).$$
    		\end{description}
    \end{cor}

We write $f_1=\Theta(f_2)$ if $f_1=\Omega(f_2)$ and $f_2=\Omega(f_1)$. Then in terms of the order of $\var(\bar{W}_\alpha)$ and $\var(W_\alpha)$, we have the following results which can be regarded as the counterparts of \cite[Theorem~1.15]{LSY19}. 

\begin{thm}\label{thmvar}
	\begin{description}
		\item{(a)} (unrestricted case) Assume that $\Xi$ satisfies the EDD and the fifth moment condition~\Ref{moment0}, and the score function $\eta$ satisfies the sixth moment condition~\Ref{thm2.1}. If $\eta$ is exponentially stabilising in Definition~\ref{defi4}, translation invariant in Definition~\ref{invar} and $\sigma^2>0$, then $\var(W_\alpha)=\Theta(\alpha)$.
		
		\item{(b)} (restricted case)  Assume that $\Xi$ satisfies the EDD and the fifth moment condition~\Ref{moment0}, and the score function $\eta$ satisfies the sixth moment condition~\Ref{thm2.1r}. If $\eta$ is exponentially stabilising in Definition~\ref{defi4r}, translation invariant in Definition~\ref{traninvres0} and $\bar{\sigma}^2>0$, then $\var(\bar{W}_\alpha)=\Theta(\alpha)$. \

		\end{description}
\end{thm}

The proofs of the main and auxiliary results are postponed to Section~\ref{Theproofs}, and we turn our attention to
the EDD point processes and applications of the main results first.

     \section{EDD point processes}\label{PPA2.0}
		
	The cornerstone model of point processes is the Poisson point process, where points behave independently in different regions, and it is obvious that a Poisson point process satisfies the EDD. There are a range of extensions of Poisson point processes to capture dependent random structures and significant development has been made in the determinantal point processes and the Gibbs point processes
\cite{SKM95,S00,Bad05,Daley03,Der19} with the connections of the two classes investigated in \cite{GY05}. Both classes have been well assessed for
statistical inferences \cite{MW04,MW07,LMR15}. In this section, we show that the Gibbs point processes with nearly finite range potentials, the determinantal point processes with fast decay kernels, the $r$-dependent point processes and the Boolean models all satisfy the EDD. 
	
\subsection{Rarified Gibbs point process}

	For ease of reading, we briefly introduce the idea of perfect simulation in \cite[Section~3]{SY13} for the Gibbs point processes with nearly finite range potentials $\Psi$. To this end, let $\Psi$ be a $[0,\infty]$ valued functional on the finite configuration space $\bm{C}_{\mathbb{R}^d,b}:=\{\xi\in\bm{C}_{\mathbb{R}^d}:\ \xi(\mathbb{R}^d)<\infty\}$ satisfying i) translation invariant: $\Psi(\mathscr{X})=\Psi(x+\mathscr{X})$ for all $x\in \mathbb{R}^d$ and $\mathscr{X}\in \bm{C}_{\mathbb{R}^d,b}$; ii) rotation invariant: $\Psi(\mathscr{X})= \Psi(\mathscr{X}')$ for all  $\mathscr{X}\in \bm{C}_{\mathbb{R}^d,b}$ and all
	rotations $\mathscr{X}'$ of $\mathscr{X}$; iii) non-decreasing: $\Psi(\mathscr{X})\le \Psi(\mathscr{X}')$ for all $\mathscr{X},\mathscr{X}'\in \bm{C}_{\mathbb{R}^d,b}$ such that $\mathscr{X}\subset \mathscr{X}'$; iv) non-degenerate:  $\Psi(\{x\})<\infty$ for all $x\in \mathbb{R}^d$. Let $D_n:=[-n,n]^d$ for $n\in \mathbb{N}:=\{1,2,\dots\}$, $\Psi_{D}(\mathscr{X}):=\Psi(\mathscr{X}\cap D)$ for $D \in \mathscr{B}(\mathbb{R}^d)$, $\mathscr{P}^{\beta\Psi}$ and $\mathscr{P}_D^{\beta\Psi}$ denote the Gibbs point process with inverse temperature $\beta>0$ and potential $\Psi$ and $\Psi_D$ respectively.  Write $\Delta(\bm{0},\mathscr{X}):=\Delta^\Psi(\bm{0},\mathscr{X}):=\Psi(\mathscr{X}\cup\{\bm{0}\})-\Psi(\mathscr{X})$, $\bm{0}\not\in \mathscr{X}$, with $\infty-\infty:=0$, and assume that  $\Delta(\bm{0},\mathscr{X})$ %, for all $n\in \mathbb{N}$, the energy function $\Delta_n(x,\mathscr{X}):=\Psi_{D_n}(\mathscr{X}\cup \{x\})-\Psi_{D_n}(\mathscr{X})$ for $\mathscr{X}\in \bm{C}_{\mathbb{R}^d,b}$, with $\infty-\infty:=0$, and $x\notin \mathscr{X}$, 
	satisfies  $$\Delta_{[r]}(\bm{0},\mathscr{X} \cap B_{r}(\bm{0}))\le \Delta(\bm{0},\mathscr{X})\le \Delta^{[r]}(\bm{0},\mathscr{X} \cap B_{r}(\bm{0}))$$ for some non-negative, translation invariant functions $\Delta_{[r]}$ and $\Delta^{[r]}$, and $\Delta_{[r]}$ and $\Delta^{[r]}$ are assumed to be respectively increasing and decreasing in $r$. We say $\beta\Psi$ has nearly finite range if there exists a decreasing continuous function $\psi^{(\beta)}: \mathbb{R}^+\rightarrow [0,1]$ satisfying that $\psi^{(\beta)}(0)=1$, $\psi^{(\beta)}(r)$ decays exponentially fast in $r$ and $$e^{-\beta \Delta_{[r]}(\bm{0},\mathscr{X}\cap B_r(\bm{0}))}-e^{-\beta \Delta^{[r]}(\bm{0},\mathscr{X}\cap B_r(\bm{0}))}\le \psi^{\beta}(r)$$ for all $r>0$ and $\mathscr{X}\in \mathbf{C}_{\mathbb{R}^d}$. \cite{SY13} established that the class of Gibbs point processes having nearly finite range $\beta\Psi$ includes 
	\begin{description}
	\item{i)} the point process with a pair potential function $\Psi(\mathscr{X})=\sum_{x\ne y}\phi(\|x-y\|)$, where $\phi:[0,\infty)\to[0,\infty)$ has a compact support or $$\phi(r)\left\{\begin{array}{ll}\le K_1\exp(-K_2r) ,\ &r\in[r_0,\infty),\\
	=\infty,\ &r\in(0,r_0),\end{array}\right.$$
	for constants $K_1,~K_2\in\mathbb{R}_+$;
	\item{ii)} the point process defined by the continuum Widom-Rowlinson model for spheres of type $A$ having centres $\mathscr{X}$ and spheres of type $B$ having centres $\mathscr{Y}$:
	$$\Psi(\mathscr{X}\cup \mathscr{Y})=\left\{\begin{array}{ll}\alpha_1\mbox{\rm card}(\mathscr{X})+\alpha_2\mbox{\rm card}(\mathscr{Y})+\alpha_3,\ &d(\mathscr{X},\mathscr{Y})>2a,\\
	\infty,\ &{\rm otherwise},\end{array}\right.$$
	where $\mbox{\rm card}(\mathscr{X})$ is the cardinality of $\mathscr{X}$, $a$ is the common radii of the spheres and $\alpha_i$'s are positive constants; 
	\item{iii)} the area interaction point process with
	$$\Psi(\mathscr{X})=\Vol\left(\cup_{x\in \mathscr{X}}(x+K)\right)+\alpha_1{\rm card}(\mathscr{X})+\alpha_2,$$
	where $\alpha_i$'s are positive constants and $K$ is a fixed compact convex set;
	\item{iv)} the hard-core process with
	$$\Psi(\mathscr{X})=\left\{\begin{array}{ll}
	\alpha_1{\rm card}(\mathscr{X})+\alpha_2,\ &\inf_{x,y\in\mathscr{X},x\ne y} |x-y|\ge r_0,\\
	\infty,\ &\inf_{x,y\in\mathscr{X},x\ne y} |x-y|<r_0,\end{array}\right.$$
	where $r_0$ and $\alpha_i$'s are positive constants. 
	\end{description}Moreover, extensions of the point process with a pair potential function have been developed in \cite{G99,RMO18} to detect multivariate interactions in spatial point patterns, and most of these extensions have nearly finite range potentials. \cite[Section~3.3]{SY13}  states that the infinite volume Gibbs point process $\mathscr{P}^{\beta \Psi}$ exists and is the thermodynamic limit of $\mathscr{P}_{D_n}^{\beta \Psi}$.
	The next lemma says that $\mathscr{P}^{\beta \Psi}$ also satisfies the EDD.
	
   \begin{lma}\label{gibbs}
   	  The Gibbs point process $\mathscr{P}^{\beta \Psi}$ with nearly finite range potential satisfies the EDD.
   \end{lma}
  
  \noindent{\it Proof.~} Using the idea of perfect simulation introduced in \cite{FFG02} and \cite[Sections~3.2 \& 3.3]{SY13}, we can construct a stationary homogeneous free birth and death process $\{\rho(t)\}_{t\in \mathbb{R}}$ such that $\rho(t)\overset{d}{=}\mathscr{P}^{\beta \Psi}$ for all $t$. For $B\in \mathscr{B}(\mathbb{R}^d)$, define the {\it ancestor clan} $\bm{A}_B^{\beta \Psi}(0)=:\bm{A}_B^{\beta \Psi}$ with respect to the process $\{\rho(t)\}_{t\in \mathbb{R}}$ as the accepted births in $\rho(0)\cap B$, their ancestors, the ancestors of their ancestors and so forth. From the construction, $\rho(0)\cap A$ and $\rho(0)\cap B$ are conditionally independent given $\bm{A}_A^{\beta \Psi}\cap \bm{A}_B^{\beta \Psi}=\emptyset$. \cite[(3.6)]{SY13} states that the ancestor clan $\bm{A}_B^{\beta \Psi}$ satisfies that for all $r\in \mathbb{R}_+$, $$\mathbb{P}\left[\mbox{diam}(\bm{A}_B^{\beta \Psi})\ge r+\mbox{diam}(B)\right]\le C(\Vol(B)\vee 1) \exp(-r/C)$$ for some positive constant $C$ depending on the distribution of the process $\mathscr{P}^{\beta \Psi}$ only. Noting that the ancestor clan $\bm{A}_B^{\beta \Psi}$ starts from the accepted births in $\rho(0)\cap B$, if $\rho(0)\cap B\ne\emptyset$, then $\bm{A}_B^{\beta \Psi}\cap B\ne\emptyset$, which ensures that if $B$ is a cube in $\mathbb{R}^d$ with centre $x$ and diagonal length $\diam(B):=2r$, then by the rotation invariance, there is a constant $C$ such that \begin{equation}\label{gibbs.1}\mathbb{P}\left[\bm{A}_B^{\beta \Psi}\nsubseteq B(x,3r+r')\right]\le C(r^d\vee 1) \exp(-r'/C)
    \end{equation} for all $r'\in \mathbb{R}_+$.
    	
  For two bounded sets $A,~B\in \mathscr{B}(\mathbb{R}^d)$ with $d(A,B)=:r_0$, without loss, we assume that $r_0\ge 1$. Let $\Xi\overset{d}{=}\mathscr{P}^{\beta \Psi}$ and $\{\tilde{\rho}(t)\}_{t\in \mathbb{R}}$ be an independent copy of $\{\rho(t)\}_{t\in \mathbb{R}}$. Since $\rho(0)\cap A$ and $\rho(0)\cap B$ are conditionally independent given $\bm{A}_A^{\beta \Psi}\cap \bm{A}_B^{\beta \Psi}=\emptyset$, we have 
    \begin{align}
    		&\beta_{A,B}\nonumber\\
		=&d_{TV}(\rho(0)\cap (A\cup B), (\rho(0)\cap A)\cup (\tilde{\rho}(0)\cap B))\nonumber
    		\\=& d_{TV}\left(\rho(0)\cap (A\cup B), (\rho(0)\cap A)\cup (\tilde{\rho}(0)\cap B)\middle| \bm{A}_A^{\beta \Psi}\cap \bm{A}_B^{\beta \Psi}=\emptyset\right)\mathbb{P}(\bm{A}_A^{\beta \Psi}\cap \bm{A}_B^{\beta \Psi}=\emptyset)\nonumber
    		\\&+d_{TV}\left(\rho(0)\cap (A\cup B), (\rho(0)\cap A)\cup (\tilde{\rho}(0)\cap B)\middle| \bm{A}_A^{\beta \Psi}\cap \bm{A}_B^{\beta \Psi}\neq\emptyset\right)\mathbb{P}(\bm{A}_A^{\beta \Psi}\cap \bm{A}_B^{\beta \Psi}\neq\emptyset)\nonumber
    		\\ \le& \mathbb{P}(\bm{A}_A^{\beta \Psi}\cap \bm{A}_B^{\beta \Psi}\neq\emptyset).\label{gibbs.2}
    \end{align} 

   Since $A$ and $B$ are bounded, $\diam(A)$ and $\diam(B)$ are finite. We can find a set of disjoint cubes $\{\mathbb{C}_{i,j}\}_{i\in\{1,2\}, 0\le j\le n_i}$ with diagonal length $\frac{r_0}{16}$ such that $A\subset \cup_{j\le n_1}\mathbb{C}_{1,j}$ and $B\subset \cup_{j\le n_2}\mathbb{C}_{2,j}$ for positive integers $n_1\le C_1(\diam(A)\vee 1)^dr_0^{-d}$ and $n_2\le C_1(\diam(B)^d\vee 1)r_0^{-d}$. Write the centre of $\mathbb{C}_{i,j}$ as $c_{i,j}$. Then for any $j_1\le n_1$, $j_2\le n_2$, $d(\mathbb{C}_{1,j_1}, \mathbb{C}_{2,j_2})\ge d(A,B)-(\diam(\mathbb{C}_{1,j_1})+\diam(\mathbb{C}_{2,j_2}))=\frac{7r_0}{8}$, it follows from \Ref{gibbs.1} with $2r=r'=\frac{r_0}{16}$ that
    \begin{align}
    	 &\mathbb{P}(\bm{A}_{\mathbb{C}_{1,j_1}}^{\beta \Psi}\cap \bm{A}_{\mathbb{C}_{2,j_2}}^{\beta \Psi}\neq\emptyset)\nonumber
    	 \\\le &\mathbb{P}\(\bm{A}_{\mathbb{C}_{1,j_1}}^{\beta \Psi}\nsubseteq B\(c_{1,j_1},\frac{5r_0}{32}\)\)+\mathbb{P}\(\bm{A}_{\mathbb{C}_{2,j_2}}^{\beta \Psi}\nsubseteq B\(c_{2,j_2},\frac{5r_0}{32}\)\)\nonumber
    	 \\ \le & C_2r_0^d \exp(-C_3r_0) \label{gibbs.3}
    \end{align} for some positive constants $C_2$ and $C_3$ independent of $j_1$ and $j_2$. Also, from the definition of the ancestor clans, if a set $B$ is covered by a class of sets $\{B_1,\dots ,B_n\}$, then $\bm{A}_{B}^{\beta \Psi}\subset \cup_{i\le n} \bm{A}_{B_i}^{\beta \Psi}$. Together with \Ref{gibbs.2} and \Ref{gibbs.3}, we have
    \begin{align}
    	\beta_{A,B}\le &\mathbb{P}(\bm{A}_A^{\beta \Psi}\cap \bm{A}_B^{\beta \Psi}\neq\emptyset)\nonumber
    	\\\le &\mathbb{P}\(\(\cup_{j_1\le n_1}\bm{A}_{\mathbb{C}_{1,j_1}}^{\beta \Psi}\)\cap \(\cup_{j_2\le n_2}\bm{A}_{\mathbb{C}_{2,j_2}}^{\beta \Psi}\)\neq\emptyset\)\nonumber
    	\\\le &\sum_{j_1\le n_1}\sum_{j_2\le n_2}\mathbb{P}\(\bm{A}_{\mathbb{C}_{1,j_1}}^{\beta \Psi}\cap \bm{A}_{\mathbb{C}_{2,j_2}}^{\beta \Psi}\neq\emptyset\)\nonumber
    	\\ \le &n_1n_2C_2r_0^{d} \exp(-C_3r_0)\le C_1^2C_2(\diam(A)^d\vee 1)(\diam(B)^d\vee 1)\exp(-C_3r_0),\nonumber
    \end{align}
    completing the proof. \qed
    
    Lemma~\ref{gibbs} ensures that Theorem~\ref{thmvar} is applicable to all geometric statistics arising from a Gibbs point process $\mathscr{P}^{\beta \Psi}$ with nearly finite range potential. In Section~\ref{Applications}, we demonstrate its use in two examples: the total edge length of a $k$-nearest neighbour graph and the total log volume in a given range of forest. 
    
	\subsection{Determinantal point process}
	
	The determinantal point processes are a broad class of point processes such that the distributions can be characterised by the determinants of given functions. More precisely, we say that $\Xi$ is a {\it determinantal point process} on space $\mathbb{R}^d$ with kernel $K$ if it is a simple point process on $\mathbb{R}^d$ with the joint intensities given by $$\rho_n(x_1,\dots,x_n)=\det\[K(x_i,x_j)\]_{1\le i,j\le n}$$ for $n\in\mathbb{N}$ and $x_1,\dots,x_n\in \mathbb{R}^d$ \cite{GY05}. These processes are widely used in random matrix theory and mathematical physics.  The next lemma says that  the determinantal point process satisfies the EDD if its kernel function decreases fast enough. 
	%Using our main result Theorem~\ref{thm2.a2}, we can show the normal approximation in  Theorem~\ref{exthm2} still holds with the vextex set being a determinantal point process with fast decay kernel under some variation condition.
	
	\begin{lma}\label{mixingdeter}
	If the kernel $K$ of the determinantal point process $\Xi$ satisfies $\|K\|_\infty:=\sup_{x,y\in\mathbb{R}^d}|K(x,y)|<\infty$ and there exist constants $C_i\in\mathbb{R}_+$, $1\le i\le 4$, such that
	$|K(x,y)|\le C_1e^{-C_2e^{C_3|x-y|}}$
	for all $x,~y\in \mathbb{R}^d$ such that $\|x-y\|\ge C_4$, then $\Xi$ is an
	EDD point process.
	\end{lma} 
	
	\noindent{\it Proof.~} We use Lemma~\ref{mixing.re} to prove the claim. Let $p_A=\Vol(A)$ and $p_B=\Vol(B)$, \cite[Theorem~4.1]{P19} states that 
    \begin{align}
    	\beta_{A,B}&\le 2p_Ap_B(1+2p_A\|K\|_\infty)(1+2p_B\|K\|_\infty)e^{2\(p_A+p_B\)\|K\|_\infty}\omega\(d(A,B)\)^2,
	\label{determinantalp1}
    \end{align} where $\omega(r):=\sup_{|x-y|\ge r}|K(x,y)|$. It is easy to see that
    $$1+2p_A\|K\|_\infty\le 3(p_A\vee 1)(\|K\|_\infty\vee 1)\mbox{ and }1+2p_B\|K\|_\infty\le 3(p_B\vee 1)(\|K\|_\infty\vee 1),$$ hence
    \begin{equation}\label{determinantalp2}
    2p_Ap_B(1+2p_A\|K\|_\infty)(1+2p_B\|K\|_\infty)
    \le 18(\|K\|_\infty^2\vee 1) (p_A^2\vee 1)(p_B^2\vee 1).
    \end{equation}
     For the remaining part, we take $C_5=\(1\vee\frac2{C_3}\)\(1\vee\ln\(\frac{4\|K\|_\infty}{C_2}\)\)$, then there exists an $r_0\in\mathbb{R}_+$ such that $s:=d(A,B)\ge C_5\ln s$ for all $s\ge r_0$. Set $\theta_3'=C_5\vee r_0$, $\theta_4'=e$ and $\Theta=p_A\vee p_B\vee e$, if $s\ge \theta_3'\ln \Theta$, then
    \begin{align}
    &e^{2\(p_A+p_B\)\|K\|_\infty}\omega\(d(A,B)\)^2\nonumber\\
    &\le C_1^2e^{4\Theta\|K\|_\infty}\exp\left\{-2C_2e^{C_3s}\right\}\nonumber\\
    &\le C_1^2e^{4\Theta\|K\|_\infty}\exp\left\{-C_2e^{C_3C_5\ln\Theta}-C_2e^{C_3C_5\ln s}\right\}\nonumber\\
   &= C_1^2e^{4\Theta\|K\|_\infty}\exp\left\{-C_2\Theta^{C_3C_5}-C_2s^{C_3C_5}\right\}\nonumber\\
    &\le C_1^2e^{4\Theta\|K\|_\infty}\exp\left\{-C_2\Theta^{1+0.5C_3C_5}-C_2s\right\}\nonumber\\
     &\le C_1^2e^{4\Theta\|K\|_\infty}\exp\left\{-C_2\Theta e^{\ln\(\frac{4\|K\|_\infty}{C_2}\) }-C_2s\right\}\nonumber\\
     &=C_1^2e^{-C_2s},\label{determinantalp3}
    \end{align}
    hence the EDD follows from combining \Ref{determinantalp1}, \Ref{determinantalp2}, \Ref{determinantalp3}
    and applying Lemma~\ref{mixing.re} with $\theta_0'=2$, $\theta_1'=18C_1^2(\|K\|_\infty^2\vee 1)$ and $\theta_2'=C_2.$     \qed
    
    Lemma~\ref{mixingdeter} ensures that Theorem~\ref{thmvar} is also applicable to all geometric statistics and the log volume arising from a determinantal point process with its kernel satisfying the conditions in Lemma~\ref{mixingdeter}. 
        
        Now, we consider two more models that satisfy the EDD so that Theorem~\ref{thmvar} can be applied to geometric statistics driven by these point processes.
        
    \subsection{$r$-dependent point process}
   A point process $\Xi$ on $\mathbb{R}^d$ is said to be {\it $r$-dependent} if for any Borel sets $A$, $B\in \mathscr{B}(\mathbb{R}^d)$ with distance $d(A,B)\ge r$, $\Xi_A$ and $\Xi_B$ are independent. The definition implies that $\beta_{A,B}=0$ for $d(A,B)\ge r$, hence the following lemma is trivial.

\begin{lma} An $r$-dependent point process satisfies the EDD.
\end{lma}

 One example of the $r$-dependent point processes is the Mat\'{e}rn hard-core process \cite[p.~298]{Daley03}. Let $\Xi':=\sum_{i\in \mathbb{N}}\delta_{X_i}$ be a homogeneous Poisson point process on $\mathbb{R}^d$. Then we can construct a hard-core Poisson process by setting $\Xi:=\sum_{i\in \mathbb{N}}\delta_{X_i}\mathbf{1}_{B(X_i,r/2)\cap \Xi'=\{X_i\}}$, then for any Borel sets $A$, $B\in \mathscr{B}(\mathbb{R}^d)$ with distance $d(A,B)\ge r$, $\sigma(\Xi_A)\subset \sigma(\Xi'_{B(A,r/2)})$ is independent of $\sigma(\Xi_B)\subset \sigma(\Xi'_{B(B,r/2)})$, where $B(A,r'):=\{x:\ d(x,A)< r'\}$ for all $A\in \mathscr{B}(\mathbb{R}^d)$ and $r'\in \mathbb{R}_+$, so $\Xi$ is $r$-dependent.

 \subsection{Boolean model}
 
 The Boolean model is a special class of the {\it germ-grain model} (see, for example, \cite{HM99}). We call a point process $\Xi=\cup_{n\in \mathbb{N}}(X_n+\Xi(n))$ a germ-grain model, where the grains $\{\Xi(n)\}_{n\in \mathbb{N}}$ are $\iid$ point processes and germs $\{X_n\}$ are independent of $\{\Xi(n)\}_{n\in \mathbb{N}}$ and form a stationary point process $\Xi':=\sum_{n\in \mathbb{N}}\delta_{X_n}$. A germ-grain model is called the {\it Boolean model} if $\Xi'$ is Poisson \cite[p.~206]{Daley03}. 

For the Boolean model, assume that $\Xi(n)$'s are bounded, that is, there exists some $r\in \mathbb{R}_+$ such that $\mathbb{P}(\Xi(n)\cap B({\bf 0},r)^c=\emptyset)=1$, then for any Borel sets $A$, $B\in \mathscr{B}(\mathbb{R}^d)$ with distance $d(A,B)\ge 4r$, $\sigma(\Xi_A)\subset \sigma\(\left\{\Xi'_{B(A,r)},\Xi(n)~\mbox{such that }X_n\in B(A,r)\right\}\)$ is independent of $\sigma(\Xi_B)$, so $\Xi$ is a $(4r)$-dependent process and, consequently, it satisfies the EDD. More generally, we have the following conclusion.
 
 \begin{lma}For the Boolean model defined above, if there exist positive constants $r_0$, $C_1$ and $C_2$ such that $\mathbb{P}(\Xi(n)\cap B({\bf 0},r)^c\neq\emptyset)\le C_1\exp(-C_2 r)$ for all $r\ge r_0$, then $\Xi$ satisfies the EDD.
 \end{lma}
 
 \noindent{\it Proof.~} Since $\beta_{A_1,A_2}$ is a non-decreasing function in the sense of inclusion (see Remark~\ref{rebeta}), without loss, we take $A_1$ and $A_2$ as two balls with centres $O_1$ and $O_2$, radii $r_1$ and $r_2$, set $R=d(A_1,A_2)$, and we assume $R\ge 4r_0$. For the points $\mathscr{X}_i$ in $A_i$ contributed by $\Xi(n)$ with its germ $X_n\in B(O_i,r_i+R/4)$ and $\Xi(n)\cap B({\bf 0},R/4)^c=\emptyset$, these points $\mathscr{X}_i$ and their dependents are in $B(O_i,r_i+R/2)$, so $\mathscr{X}_1$ and $\mathscr{X}_2$ are independent. In other words, the contribution of dependent points is from those violating these conditions. By abuse of notation, we write $\Xi(x)$ as the grain of the germ at $x\in\Xi'$, $V(r)=\frac{\pi^{d/2}}{\Gamma(1+d/2)}r^d$ as the volume of the ball with radius $r$, and add up the probabilities of all possible cases leading to dependent points in $A_1$ and $A_2$ to get
 \begin{align*}
 \beta_{A_1,A_2}\le &\sum_{i=1}^2\int_{B(O_i,r_i+R/4)}\mathbb{P}\(\Xi(x)\cap B({\bf 0},R/4)^c\ne\emptyset\)\mathbb{E}\Xi'(dx)\\
 &+\sum_{i=1}^2\int_{B(O_i,r_i+R/4)^c} \mathbb{P}\(\Xi(x)\cap B({\bf 0},\|x-O_i\|-r_i)^c\ne\emptyset\)\mathbb{E}\Xi'(dx)\\
\le &\sum_{i=1}^2 C_1\lambda e^{-C_2R/4}V(r_i+R/4)+\sum_{i=1}^2\int_{r_i+R/4}^\infty C_1 \lambda e^{-C_2(r-r_i)}dV(r)\\
 \le&C_3e^{-C_2R/4}(R\vee r_1\vee r_2\vee 1)^d,
 \end{align*}
 for $C_3\in\mathbb{R}_+,$ where $\lambda=\mathbb{E}\Xi'(d{\bf 0})/d{\bf 0}$.
 Choose $R_0\ge 1\vee (4r_0)$ such that $\frac{R}{\ln R}\ge \frac{8d}{C_2}$ for all $R\ge R_0$, set $\theta_4=e^{R_0}$, $\theta_3=1$, then for $R\ge \ln((2r_1)\vee (2r_2)\vee \theta_4)$, we have
 $$\beta_{A_1,A_2}\le C_3((2r_1)^d\vee 1)((2r_2)^d\vee 1)e^{-C_2R/8},$$
 which implies \Ref{mixing} with $\theta_0=d$, $\theta_1=C_3$, $\theta_2=C_2/8$.
 \qed

	\section{Applications}\label{Applications}
	The asymptotic behaviour of geometric functionals has been of considerable interest in the last three decades, and our main normal approximation results can be applied to a large class of geometric functionals, including $k$-nearest neighbour graph, Voronoi graph, sphere of influence graph, Delaunay triangulation, Gabriel graph and relative neighbourhood graph \cite{Devroye88, T82} with vertices driven by a point process satisfying the EDD. 
	
The limit theory of geometric functionals with determinantal point process input or Gibbs point process input is investigated in \cite{BYY19}. Error bounds of a normal approximation in terms of the Kolmogorov distance for the geometric functions with Gibbsian input were derived in \cite{XY15}. 
To illustrate the use of the main results in Section~\ref{Generalresults}, we bound the errors of a normal approximation to the total edge length in the $k$-nearest neighbour graph with vertices forming a rarified Gibbs point process or a determinantal point process with fast decay kernel. We also bound the error of a normal approximation to the total log volume in a given range of forest with trees following a marked Gibbs point process.
	
 \subsection{The total edge length of $k$-nearest neighbour graphs}%\label{sec4.1}
    
    The $k$-nearest neighbour graph $NG(\mathscr{X})$ with respect to a configuration $\mathscr{X} \in \bm{C}_{\mathbb{R}^d}$ is a graph with vertices $\mathscr{X}$  and edges $\{x,y\}$ such that $y$ is one of the $k$ points nearest to $x$ or $x$ is one of the $k$ points nearest to $y$ in $\mathscr{X}$. A variant $NG'\(\mathscr{X}\)$ of the $NG(\mathscr{X})$ can be constructed by inserting directed edges $x\rightarrow y$ if $y$ is one of the $k$ nearest neighbours of $x$ instead of the undirected edges in $NG(\mathscr{X})$. As in \cite{SY13}, we take the score function $\eta(x,\mathscr{X},\Gamma_{\alpha})$ (resp. $\eta'(x,\mathscr{X},\Gamma_{\alpha})$) as one half the sum of the edge lengths of edges in $NG\left(\Gamma_{\alpha}\cap\left(\mathscr{X}\cup\{x\}\right)\right)$ (resp. $NG'\left(\Gamma_{\alpha}\cap\left(\mathscr{X}\cup\{x\}\right)\right)$) which are
incident to $x$, and set
    \begin{equation}\bar{W}_\alpha=\sum_{{x\in\Xi_{\Gamma_\alpha}}}\eta(x, \Xi,\Gamma_\alpha)\mbox{ and } \bar{W}_\alpha'=\sum_{{x\in\Xi_{\Gamma_\alpha}}}\eta'(x, \Xi,\Gamma_\alpha).\label{defexam1}\end{equation}
  We now state the error bounds for a normal approximation of the total edge lengths $\bar{W}_\alpha$ and $\bar{W}_\alpha'$ of $NG\(\Xi_{\Gamma_{\alpha}}\)$ if $\Xi$ follows $\mathscr{P}^{\beta \Psi}$ or a determinantal point process with fast decay dependence.
	 
	\begin{thm}\label{exthm2}
		\begin{description}
		\item{(a)} If $\Xi$ is an infinite range Gibbs point process with nearly finite range potential, then  
		$$d_W\(\frac{\bar{W}_{\alpha}-\mean \bar{W}_{\alpha}}{\sqrt{\var(\bar{W}_{\alpha})}},Z\)\le O\(\alpha^{-\frac{1}{2}}\ln(\alpha)^{5d}\).$$
		The statement holds if $\bar{W}_\alpha$ is replaced by $\bar{W}_\alpha'$ in \Ref{defexam1}.
		\item {(b)} 		If $\Xi$ is a determinantal point process with continuous kernel $K$ satisfying the conditions in Lemma~\ref{mixingdeter}, and the total edge length $\bar{W}_\alpha$ satisfies $\var(\bar{W}_\alpha)= \Omega(\alpha^{\nu})$ for some $\nu>\frac{2}{3}$, then
		$$d_W\(\frac{\bar{W}_{\alpha}-\mean \bar{W}_{\alpha}}{\sqrt{\var(\bar{W}_{\alpha})}},Z\)\le O\(\alpha^{-\frac{3}{2}\nu+1}\ln(\alpha)^{5d}\).$$
		The statement holds true if $\bar{W}_\alpha$ is replaced by $\bar{W}_\alpha'$ in \Ref{defexam1}.
		\end{description}
	\end{thm}
	 
	 \begin{re} \cite{XY15} proved that a normal approximation error of $W_\alpha:=\sum_{{x\in\Xi_{\Gamma_\alpha}}}\eta(x, \Xi,\Gamma_\infty)$ in terms of the Kolmogorov distance can be bounded above by
	 $O\(\alpha^{-\frac{1}{2}}\ln(\alpha)^{2d}\)$, which is slightly better than the error bound for $\bar{W}_\alpha$ in terms of the Wasserstein distance in Theorem~\ref{exthm2}~(a).
	 \end{re}
	 
	 \noindent{\it Proof of Theorem~\ref{exthm2}.~} We only show the claims for the undirected case, and the directed case can be handled using the same idea. In this case, the total edge length $\bar{W}_{\alpha}$ can be represented as 
	  $$\eta\(x,\mathscr{X},\Gamma_\alpha\):=\frac{1}{2}\sum_{y\in \mathscr{X}_{\Gamma_\alpha}}\|y-x\|\mathbf{1}_{\{(x,y)\in NG(\mathscr{X}_{\Gamma_\alpha})\}}.$$
	   The score function $\eta$ is translation invariant according to the construction and the EDD is ensured by Lemma~\ref{gibbs} for (a) and Lemma~\ref{mixingdeter}  for (b). To apply Theorem~\ref{thm2.a2}, we need to check the stabilising condition as in Definition~\ref{defi4r}, the moment conditions \Ref{moment0}, \Ref{thm2.1r} and the order of $\Var(\bar{W}_{\alpha}).$
	 
	 (a) According to Lemma~$3.3$ in \cite{SY13}, the Gibbs point process $\Xi$ is {\it Poisson-like}, which means that $\Xi$ is stochastically dominated by a Poisson point process on $\mathbb{R}^d$ with intensity $\lambda'>0$, and there exist strictly positive constants $C:=C(\lambda')$ and $r_1$ such that for all $r\ge r_1$, $x\in \mathbb{R}^d$ and $\mathscr{X}\in \bm{C}_{\mathbb{R}^d\backslash B(x,r)}$, the conditional probability that $B(x,r)$ is not hitted by $\Xi$ given $\Xi_{B(x,r)^c}=\mathscr{X}$ satisfies that \begin{equation}\mathbb{P}\left[\Xi_{B(x,r)}=\emptyset\middle|\Xi_{B(x,r)^c}=\mathscr{X}\right]\le e^{-Cr^d}.\label{nearest0}\end{equation}
	 Consequently, the fifth moment condition~\Ref{moment0} of $\Xi$ is ensured by the Poisson-like property and the moment property of the Poisson point process. 
	 
	To examine the remaining conditions, for simplicity, we take $d=2$ and follow the proof of Theorem~3.1 in \cite{CX20} using the idea initiated in \cite{PY01} to achieve the purpose. For completeness, we recap the main steps in \cite{CX20}. For $x\in\Gamma_\alpha$ and $t>0$, we carve the disk with centre $x$ and radius $t$ into six disjoint circular sectors $T_j(t)$, $1\le j\le 6$, of the same size with $x$ as the centre and $\frac{\pi}{3}$ as their central angle. The sectors are rotated around $x$ such that all straight edges of the sectors have at least the minimal angle $\pi/12$ with respect to the edges of $\Gamma_\alpha$. Let $T_{j}(\infty)=\cup_{t>0}T_j(t)$ for $1\le j\le 6$ and define
	$$t_{x,\alpha}(\Xi_x)=\inf\{t:\ \mbox{\rm card}(T_j(t)\cap \Gamma_\alpha\cap \Xi_x)\ge k+1\mbox { or }T_j(t)\cap \Gamma_\alpha=T_j(\infty)\cap \Gamma_\alpha,\ 1\le j\le 6\}$$
	and $\bar{R}\(x,\alpha\)=3t_{x,\alpha}(\Xi_x+\delta_x)$. It was demonstrated in \cite{CX20} that $\bar{R}$ is a radius of stabilisation. For the tail distribution of $\bar{R}$, let $A_t$ be an obtuse triangle with the longest side length $t$ and two angles $\pi/12$ and $\pi/3$, define $\tau:=\inf\{t:\ \mbox{\rm card}(\Xi_x\cap A_t)\ge k+1\}$, then \cite{CX20} established that $\mathbb{P}\(\bar{R}\(x,\alpha\)>t\)\le 6\mathbb{P}\(\tau>t/3\)$. We can find a constant $C_1\in\mathbb{R}_+$ such that there are $k+1$ disjoint disks $\{B_1,\dots,B_{k+1}\}$ of radius $C_1t$ such that $\cup_{i=1}^{k+1}B_i\subset A_{t/3}$. Using Poisson-like property \Ref{nearest0}, we obtain
	 \begin{align}&\mathbb{P}\(\bar{R}\(x,\alpha\)>t\)\le 6\mathbb{P}\(\tau>t/3\)\nonumber\\
	 &\le 6\mathbb{P}\(\mbox{\rm card}(\Xi_x\cap A_{t/3})\le k\)
	 \le 6\mathbb{P}\(\cup_{i=1}^{k+1}\{\Xi_x\cap B_i =\emptyset\}\)\nonumber\\
	 &\le 6(k+1)e^{-C\(C_1t\)^2},\label{nearest2}
	 \end{align}
	which ensures the exponential stabilisation in Definition~\ref{defi4r}. {For the moment condition~\Ref{thm2.1r}, we again make use of the proof of Theorem~3.1 in \cite{CX20} that \begin{equation}\eta\(x,(\Xi_x)_{\Gamma_{\alpha}}+\delta_x\)\le 3.5kt_{x,\alpha}(\Xi_x+\delta_x).
			\label{nearest3}
	\end{equation} Since $\bar{R}\(x,\alpha\)=3t_{x,\alpha}(\Xi_x)$, \Ref{nearest2} implies
\begin{align*}\sup_{\alpha\in \mathbb{R}_+}\sup_{x\in \Gamma_\alpha}\mathbb{P}\(t_{x,\alpha}(\Xi_x+\delta_x)>t\)\le 6(k+1)e^{-C\(3C_1t\)^2}
\end{align*}
for all $t>0$.  This ensures $\sup_{\alpha\in\mathbb{R}_+}\sup_{x\in\Gamma_\alpha}\mathbb{E}\(t_{x,\alpha}(\Xi_x+\delta_x)\)^6<\infty$ and the moment condition~\Ref{thm2.1r} is an immediate consequence of \Ref{nearest3}.
Finally, we establish \Ref{non-sinr}. To this end, define  $$\bar{\eta}\(x,\mathscr{X}\):=\frac{1}{2}\sum_{y\in \mathscr{X}}\|y-x\|\mathbf{1}_{\{(x,y)\in NG(\mathscr{X})\}},$$
and  $\bar{W}_{\infty,\alpha}=\sum_{{x\in\Xi_{\Gamma_\alpha}}}\bar{\eta}(x, \Xi)$,
then we can apply \cite[Theorem~1.1]{XY15} to obtain \Ref{non-sinr}.  The proof of (a) is completed by applying Corollary~\ref{thm2a}~(ii).

(b) Since the kernel $K$ is continuous and fast-decreasing, it follows from Section~2.2.2 and Section~2.1, Remark~(i) of \cite{BYY19} that $\mathbb{E}\(\Xi(B)^k\)$ is finite for all bounded $B\in \mathscr{B}(\mathbb{R}^d)$ and all $k\in \mathbb{N}$, which ensures the fifth moment condition~\Ref{moment0} of $\Xi$.
    
To apply Theorem~\ref{thm2.a2}, as the order of $\Var(\bar{W}_{\alpha})$ is assumed, it remains to check the stabilising condition as in Definition~\ref{defi4r} and the moment condition \Ref{thm2.1r}. For simplicity, we again take $d=2$. Following the same argument as that for \Ref{nearest2} and applying Lemma~5.6 of \cite[Supplement]{BYY19}, we obtain 
 \begin{align}\mathbb{P}\(\bar{R}\(x,\alpha\)>t\)\le 6\mathbb{P}\(\cup_{i=1}^{k+1}\{\Xi_x\cap B_i =\emptyset\}\)\le 6(k+1)e^{1/8-K({\bf 0},{\bf 0})\pi \(C_1t\)^2/8},\label{nearest4}
	 \end{align}
	 which implies the exponential stabilisation as in Definition~\ref{defi4r}.  For the moment condition \Ref{thm2.1r}, we again use the relationship $\bar{R}\(x,\alpha\)=3t_{x,\alpha}(\Xi_x)$ and \Ref{nearest4} to get 
	 $$\sup_{\alpha\in \mathbb{R}_+}\sup_{x\in \Gamma_\alpha}\mathbb{P}\(t_{x,\alpha}(\Xi_x+\delta_x)>t\)\le 6(k+1)e^{1/8-K({\bf 0},{\bf 0})\pi \(3C_1t\)^2/8}$$ for all  $t>0$. The tail behaviour of $t_{x,\alpha}(\Xi_x+\delta_x)$ and \Ref{nearest3} ensure that
	 $$\sup_{\alpha\in \mathbb{R}_+}\sup_{x\in \Gamma_\alpha}\mathbb{E}\(\eta\(x,(\Xi_x)_{\Gamma_{\alpha}}+\delta_x\)^6\)\le (3.5k)^6\sup_{\alpha\in \mathbb{R}_+}\sup_{x\in \Gamma_\alpha}\mathbb{E}\(t_{x,\alpha}(\Xi_x+\delta_x)^6\)<\infty.$$} The proof of (b) is completed by applying  Theorem~\ref{thm2.a2}~(ii).\qed
	 
	 	\subsection{The log volume of a forest with Gibbs point process tree locations}

Marks play an important role when it is necessary to classify the points. For example, in insurance, marks may be introduced to represent the types of claims \cite{ZS21}; in thinning~\cite[p.~32]{Daley03}, marks may be used to stand for the points retained and discarded. In this subsection, we consider the total log volume in a random forest \cite[Section~$3.3$]{CX20}, where marks are used to label the species of the trees. The estimation of the total log volume in a given range is of great interest in forest science and forest management \cite{C80,Li15}. When modelling the natural forest, it is reasonable to assume the locations of trees form a Gibbs point process $\overline{\Xi}$, such as a Poisson point process or a hard-core process. As the contribution of the log volume from different  species of trees varies, we use marks to classify the species. That is, for $x\in \overline{\Xi}$, let  $M_x\in{T}:=\{1,\dots,n\}$ to be the species of the tree at position $x$. We can assume that the marks are independent of other marks and the locations $\overline{\Xi}$. Then $\Xi:=\sum_{x\in \overline{\Xi}}\delta_{(x,M_x)}$ forms a marked Gibbs point process recording the locations and species of trees in a forest. We can model the timber volume of the tree at location $x$ by a function of the location, the species of the tree and the configuration of trees in a finite range around $x$, adjusted by a quantity $\e_x$ due to other unspecified factors. Formally speaking, the timber volume of a tree at location $x$ can be denoted by $\(\eta(\left(x,m\right), \Xi_{\Gamma_\alpha},\Gamma_\alpha)+\e_x\)\vee 0$, where $\eta$ is a non-negative bounded score function such that $$\eta(\left(x,m\right), \Xi_{\Gamma_\alpha},\Gamma_\alpha)=\eta(\left(x,m\right), \Xi_{\Gamma_\alpha\cap B(x,r)},\Gamma_\alpha)$$ for some positive constant $r$ and $$\eta(\left(x,m\right), \Xi_{B(x,r)},\Gamma_{\alpha_1})=\eta(\left(x,m\right), \Xi_{B(x,r)},\Gamma_{\alpha_2})$$ for all $\alpha_1$ and $\alpha_2$ with $B(x,r)\subset \Gamma_{\alpha_1\wedge\alpha_2}$. Then we have the following result analogous to \cite[Theorem~$3.3$]{CX20}.
	\begin{thm} If $\Xi$ is an infinite range Gibbs point process with nearly finite range potential, $\e_{x}$'s are $\iid$ random variables with finite sixth moment and the positive part $\e_x^+:=\e_x\vee 0$ is non-degenerate (i.e., $\var(\e_x^+)>0$), and $\e_{x}$'s are independent of $\Xi$, then the log volume in the range $\Gamma_\alpha$ is
		$${\bar{W}}_{\alpha}:=\sum_{x\in \overline{\Xi}_{\Gamma_\alpha}}\(\eta(\left(x,m\right), \Xi_{\Gamma_\alpha},\Gamma_\alpha)+\e_x\)\vee 0$$
	 and it satisfies
		$$d_W\(\frac{\bar{W}_\alpha-\mean \bar{W}_{\alpha}}{\sqrt{\var(\bar{W}_{\alpha})}},Z\)\le O\(\alpha^{-\frac{1}{2}}\).$$
\end{thm}
\noindent{\it Proof.~} The proof is adapted from  \cite[Theorem~$3.3$]{CX20}. We can construct a new marked Gibbs point process $\Xi':=\sum_{x\in \overline{\Xi}}\delta_{(x, (M_x, \e_x))}$ by replacing the marks $\{M_x\}_{x\in \bar{\Xi}}$ of $\Xi$ by $\iid$ marks $\{(M_x,\e_x)\}_{x\in \bar{\Xi}}$ on the space $(T\times \mathbb{R}, \mathscr{T}\times \mathscr{B}(\mathbb{R}))$ independent of the ground process $\overline{\Xi}'=\overline{\Xi}$. The fifth moment condition~\Ref{moment0} of $\Xi'$ follows from the Poisson-like property, as shown in the proof of Theorem~\ref{exthm2}~(a), and the EDD is ensured by Lemma~\ref{gibbs}. As $\bar{W}_\alpha$ can be represented as the sum of score function
 \begin{align*}\eta'((x,{(m,\e_x)}),\Xi',\Gamma_\alpha):=& \eta'((x,(m,\e_x)),\Xi'_{\Gamma_\alpha},\Gamma_\alpha)\\
 	:=&\[\(\eta(\left(x,m\right), \Xi_{\Gamma_\alpha},\Gamma_\alpha)+\e_x\)\vee 0\]\mathbf{1}_{(x,(m,\e_x))\in \Xi'_{\Gamma_\alpha}},
\end{align*}  the translation invariant property and the range-bound property are direct results of the construction, the sixth moment condition~\Ref{thm2.1r} is guaranteed by the boundedness of $\eta$, the moment condition of $\e_x$'s and the Minkowski inequality. Now, we can apply \cite[Theorem~1.1]{XY15} again to get $\Var(\bar{W}_{\alpha})=\Omega(\alpha)$ and the proof is then completed by applying Theorem~\ref{thm2.a2}~(i). \qed

	\begin{re} {\rm As discussed in \cite[Remark $3.2$]{CX20}, if the timber volume is determined by its $k$-nearest neighbouring trees, we can adjust the above proof to show the distribution of the log volume $\bar{W}_{\alpha}$ satisfies $d_W\(\frac{\bar{W}_{\alpha}-\mean \bar{W}_{\alpha}}{\sqrt{\var(\bar{W}_{\alpha})}},Z\)\le O\(\alpha^{-\frac{1}{2}}\ln(\alpha)^{5d}\).$}  \end{re}
    
\section{The proofs of the auxiliary and main results}\label{Theproofs}

Recalling the shift operator defined in Section~\ref{Generalresults}, we can write $g(\mathscr{X}^x):=\eta(\left({\bf 0},m\right), \mathscr{X}^x)$ $=\eta(\left(x,m\right), \mathscr{X})$ (resp. $g_\alpha(x,\mathscr{X}):= \eta\(\(x,m\),\mathscr{X},\Gamma_\alpha\)$) for all configurations $\mathscr{X}$ with $(x,m)\in \mathscr{X}$ and $\alpha>0$ so that notations can be simplified, e.g.,
\begin{eqnarray*}
		W_\alpha&=&\sum_{(x,m)\in\Xi_{\Gamma_\alpha}}\eta(\left(x,m\right), \Xi)=\int_{\Gamma_\alpha} g(\Xi^x)\overline{\Xi}(dx)=\sum_{x\in \overline{\Xi}_{\Gamma_\alpha}}g(\Xi^x),\\
		\bar{W}_\alpha&=&\sum_{(x,m)\in\Xi_{\Gamma_\alpha}}\eta(\left(x,m\right), \Xi,\Gamma_\alpha)=\int_{\Gamma_\alpha}{g_\alpha(x, {\Xi})\overline{\Xi}(dx)}=\sum_{x\in \overline{\Xi}_{\Gamma_\alpha}}g_\alpha(x, {\Xi}),
\end{eqnarray*}
where $\overline{\Xi}$ is the projection of $\Xi$ on $\mathbb{R}^d$.

We now proceed to establish a few lemmas needed in the proofs. The following lemma bounds the difference between a normal distribution and the standard normal distribution under the Wasserstein distance, and it can be verified directly (see also \cite[Lemma~2.4]{CM10}).
\begin{lma}\label{lma10}
	Let $F_{\mu,\sigma}$ be the distribution of $N(\mu,\sigma^2)$, the normal distribution with mean $\mu$ and variance $\sigma^2$, and $\Phi=F_{0,1}$, then 
	$$d_W(F_{\mu,\sigma},\Phi)\le |\mu|+\frac{2}{\sqrt{2\pi}}|\sigma-1|.$$
\end{lma}

The following lemma says that the cost of throwing away the terms with large radii of stabilisation is negligible under stabilising conditions. For convenience, we define $W_{\alpha,r}:=\sum_{(x,m){\in \Xi_{\Gamma_\alpha}}}\eta(\left(x,m\right), \Xi)\mathbf{1}_{R(x)\le r}$, $\bar{W}_{\alpha,r}:={\sum_{(x,m)\in \Xi_{\Gamma_\alpha}}\eta(\left(x,m\right), \Xi,\Gamma_\alpha)\mathbf{1}_{\bar{R}(x,\alpha)\le r}}$, which means that we through away the terms with stabilisation radii greater than $r$ from $W_\alpha$ and $\bar{W}_\alpha$.

\begin{lma}\label{lma105} 
	\begin{description}
		\item{(a)} (unrestricted case) If the score function is exponentially stabilising in Definition~\ref{defi4}, then we have
		$$d_{TV}(W_\alpha, W_{\alpha,r})\le C_1\alpha e^{-C_2r}$$
		for some positive constants $C_1$, $C_2$.
		
		\item{(b)} (restricted case) If the score function is exponentially stabilising in Definition~\ref{defi4r}, then we have
		$$d_{TV}(\bar{W}_\alpha, \bar{W}_{\alpha,r})\le C_1\alpha e^{-C_2r}$$
		for some positive constants $C_1$, $C_2$.
	\end{description}
\end{lma}
\noindent{\it Proof.} We first prove (b). Recall that $M_x\sim\mathscr{L}_T$ is the mark of the point $x\in\overline{\Xi}$, and it is independent of $\Xi_x$. From the construction of $\bar{W}_\alpha$ and $\bar{W}_{\alpha,r}$, we can see that the event $\{\bar{W}_\alpha\neq \bar{W}_{\alpha,r}\}\subset\{\mbox{at least one }x\in \overline{\Xi}\cap\Gamma_\alpha~\mbox{with }\bar{R}(x, \alpha)>r\}$, so from \Ref{palm4}, we have 
\begin{align}
	d_{TV}(\bar{W}_\alpha, \bar{W}_{\alpha,r})\le &\mathbb{P}\(\{\bar{W}_\alpha\neq \bar{W}_{\alpha,r}\}\)\nonumber
	\\\le & \mathbb{P}\(\{\mbox{at least one }x\in \overline{\Xi}\cap\Gamma_\alpha~{\mbox{such that }}\bar{R}(x, \alpha)>r\}\)\nonumber
	\\ \le& \mathbb{E}\int_{\Gamma_\alpha}\mathbf{1}_{\bar{R}(x,\alpha)>r} \overline{\Xi}(dx)\nonumber
%	\\ =&\int_{\Gamma_\alpha}\mathbb{E}\(\mathbf{1}_{\bar{R}(x,M_x,\alpha, \Xi_x+\delta_{(x, M_x)})>r}\)\lambda dx\nonumber
	\\ =&\int_{\Gamma_\alpha}\mathbb{P}\(\bar{R}((x,M_x),\alpha, \Xi_x+\delta_{(x, M_x)})>r\) \lambda dx\nonumber
	\\ \le&\alpha\lambda\bar{\tau}(r).\label{lma105.1}
\end{align} This, together with the stabilisation condition in Definition~\ref{defi4r}, gives the claim in (b).

The claim (a) can be proved by replacing corresponding counterparts $\bar{W}_\alpha$ by $W_\alpha$; $\bar{W}_{\alpha,r}$ by $W_{\alpha,r}$; $\bar{R}(x,\alpha)$ by $R(x)$; $\bar{R}((x,M_x),\alpha,\Xi_x+\delta_{(x, M_x)})$ by $R((x,M_x),\Xi_x+\delta_{(x, M_x)})$ and $\bar{\tau}$ by $\tau.$\qed

The moments of $W_{\alpha,r}$ and $W_\alpha$ (resp. $\bar{W}_{\alpha,r}$ and $\bar{W}_\alpha$) can be established using the moment conditions required. To begin with, we first show a statement about the moments of $\overline{\Xi}(\Gamma_\alpha)$. Let $\|X\|_p:=\mathbb{E}\(|X|^p\)^{\frac{1}{p}}$ be the $L_p$ norm of $X$ provided it is finite.

\begin{lma}\label{lma1} For $k\in\mathbb{N}$, if the marked point process $\Xi\sim\mathscr{P}$ satisfies the $k$th moment condition \Ref{moment0}, then $\mathbb{E}\(\overline{\Xi}(\Gamma_\alpha)^k\)\le O(\alpha^k)$ for $\alpha>0$.
\end{lma}
\noindent{\it Proof.} Since $\overline{\Xi}(B)$ is non-decreasing in $B$ in the sense of inclusion and it is also stationary,  the condition \Ref{moment0} is equivalent to that there exists an $\alpha_0>0$ such that $\mathbb{E}\(\overline{\Xi}(\Gamma_{\alpha_0})^k\)=\mathbb{E}\(\overline{\Xi}(\Gamma_{\alpha_0}+x)^k\)=:C< \infty$ for all $x\in \mathbb{R}^d$.

We can find a cover of $\Gamma_\alpha$ of the form $\{\Gamma_{\alpha_0}+x_i\}_{i\le n_\alpha}$ with $n_\alpha=\left\lceil\(\frac{\alpha}{\alpha_0}\)^\frac{1}{d}\right\rceil^d=O(\alpha)$, it follows from Minkowski's inequality that
\begin{align*}
	\mathbb{E}\(\overline{\Xi}(\Gamma_\alpha)^k\)=\|\overline{\Xi}(\Gamma_\alpha)\|_k^k\le \(\sum_{i=1}^{n_\alpha}\|\overline{\Xi}(\Gamma_{\alpha_0}+x_i)\|_k\)^k=n_\alpha^k C^k=O(\alpha^k),
\end{align*}
as claimed. \qed

\begin{re}\label{remoment} From the proof above, we can see that for arbitrary $A\in \mathscr{B}\(\mathbb{R}^d\)$, if $A$ can be covered by $\{\Gamma_{\alpha_0}+x_i\}_{i\le n_A}$, then $\mathbb{E}\(\overline{\Xi}(A)^k\)\le n_A^k C$ for some constant $C$.
\end{re}

\begin{lma}\label{lma11} 
	\begin{description}
		\item{(a)} (unrestricted case) If $\Xi$ satisfies the $(2n-1)$th moment condition~\Ref{moment0} and the score function $\eta$ satisfies the $(2n)$th moment condition~\Ref{thm2.1}, then
		$$\mathbb{E}\(|W_\alpha|^k\)\vee \mathbb{E}\(|W_{\alpha,r}|^k\)\le C\alpha^{2k-1}$$ for some positive constant $C$ for all integers $1\le k\le n$.
		\item{(b)} (restricted case) If $\Xi$ satisfies the $(2n-1)$th moment condition~\Ref{moment0} and the score function $\eta$ satisfies the $(2n)$th moment condition~\Ref{thm2.1r}, then
		$$\mathbb{E}\(|\bar{W}_\alpha|^k\)\vee \mathbb{E}\(|\bar{W}_{\alpha,r}|^k\)\le C\alpha^{2k-1}$$ for some positive constant $C$ for all integers $1\le k\le n$.
	\end{description}
\end{lma}
\noindent{\it Proof.} We start with the restricted case, and the unrestricted case follows in the same way. For $k\in \mathbb{N}$,
\begin{align}
	\mathbb{E}\(|\bar{W}_\alpha|^k\)=&\mathbb{E}\(\left|\int_{\Gamma_\alpha}g_{\alpha}(x,\Xi)\overline{\Xi}(dx)\right|^k\)\nonumber
	\\\le&\mathbb{E}\(\(\int_{\Gamma_\alpha}\left|g_{\alpha}(x,\Xi)\right|\overline{\Xi}(dx)\)^k\)\nonumber
	\\ = &\mathbb{E}\underbrace{\int_{\Gamma_\alpha}\dots\int_{\Gamma_\alpha}}_{k~\mbox{\scriptsize{of them}}}\left|g_{\alpha}(x_1,\Xi)\dots g_{\alpha}(x_k,\Xi)\right|\overline{\Xi}(dx_1)\dots \overline{\Xi}(dx_k)\nonumber
	\\ \le &\frac1k\mathbb{E}\underbrace{\int_{\Gamma_\alpha}\dots\int_{\Gamma_\alpha}}_{k~\mbox{\scriptsize{of them}}}\(\left|g_{\alpha}(x_1,\Xi)\right|^k+\dots +\left|g_{\alpha}(x_k,\Xi)\right|^k\)d\overline{\Xi}(dx_1)\dots \overline{\Xi}(dx_k)\nonumber
	\\=&\mathbb{E}\int_{\Gamma_\alpha}\overline{\Xi}(\Gamma_\alpha)^{k-1}\left|g_{\alpha}(x,\Xi)\right|^k\overline{\Xi}(dx)\nonumber
	\\ \le &\frac{1}{2}\mathbb{E}\int_{\Gamma_\alpha}\(\overline{\Xi}(\Gamma_\alpha)^{2k-2}+\left|g_{\alpha}(x,\Xi)\right|^{2k}\)\overline{\Xi}(dx)\nonumber
	\\=&\frac{1}{2}\[\mathbb{E}\(\overline{\Xi}(\Gamma_\alpha)^{2k-1}\)+\mathbb{E}\(\int_{\Gamma_\alpha}|g_\alpha(x,\Xi_x+\delta_{(x,M_x)})|^{2k}\lambda dx\)\]\nonumber
	\\\le & \frac{1}{2}C_1\alpha^{2k-1}+\frac{1}{2}C_2\alpha, \label{lma11.1}
\end{align}
for some positive constants $C_1$ and $C_2$, the second and third inequalities follow from the fact that $\Pi_{i=1}^jy_i\le \frac1j \sum_{i=1}^j y_i^j$ for all $y_i\ge 0$, $1\le i\le j$, the last equality follows from \Ref{palm4} and the last inequality follows from Lemma~\ref{lma1}. Then we can find a common positive constant $C$ such that $\mathbb{E}\(|\bar{W}_\alpha|^k\)\le C\alpha^{2k-1}$ for all integers $k\le n$. The claim $\mathbb{E}\(|\bar{W}_{\alpha,r}|^k\)\le C\alpha^{2k-1}$ can be proved by following exactly the same steps but replacing $g_{\alpha}(x,\Xi)$ with $g_{\alpha}(x,\Xi)\mathbf{1}_{\bar{R}(x,\alpha)\le r}$.
    
The statement for the unrestricted case is also true, which can be proved by replacing the corresponding counterparts $\bar{W}_{\alpha}$ by $W_{\alpha}$, $\bar{W}_{\alpha,r}$ by $W_{\alpha,r}$, and $g_\alpha$ by $g$. \qed

\vskip10pt
\begin{re}\label{moment2} {\rm~The proof of Lemma~\ref{lma11} does not depend on the shape of $\Gamma_\alpha$, so the claims still hold if we replace $\Gamma_\alpha$ with a set $A\in\mathscr{B}(\mathbb{R})$ such that $A$ satisfies the assumption in Remark~\ref{remoment} with $n_A\le O(\Vol(A))$}.
\end{re}
\begin{re} \label{moment1} {\rm~Using the same idea as in the proof of \Ref{lma11.1}, for each $1\le i\le k$, by taking the range of $x_i$ in $A_i$ satisfying the condition in Remark~\ref{moment2} instead of $\Gamma_\alpha$, we have, for all integers $1\le k\le n$,
\begin{align*}\mathbb{E}\int_{A_1}\dots\int_{A_k}\left|g_{\alpha}(x_1,\Xi)\dots g_{\alpha}(x_k,\Xi)\right|\overline{\Xi}(dx_1)\dots \overline{\Xi}(dx_k)\le C\max_{1\le j\le k}\Vol (A_j)^{2k-1} \end{align*} for some constant $C$.} 
\end{re}

With these preparations, we are ready to bound the differences $\left|\var\(W_{\alpha}\)-\var\(W_{\alpha,r}\)\right|$ and $\left|\var\(\bar{W}_{\alpha}\)-\var\(\bar{W}_{\alpha,r}\)\right|$.

\begin{lma}\label{lma12}
	\begin{description}
		\item{(a)} (unrestricted case) Assume that $\Xi$ satisfies the fifth moment condition~\Ref{moment0} and the score function $\eta$ satisfies the sixth moment condition~\Ref{thm2.1}. If $\eta$ is exponentially stabilising in Definition~\ref{defi4}, then there exist positive constants $\alpha_0$ and $C$ such that 
		\begin{equation}\left|\var\(W_{\alpha}\)-\var\(W_{\alpha,r}\)\right|\le \frac{1}{\alpha}\label{lma12s1}\end{equation} for all $\alpha\ge\alpha_0$ and $r\ge C\ln (\alpha)$.
		\item{(b)} (restricted case) Assume that $\Xi$ satisfies the fifth moment condition~\Ref{moment0} and the score function $\eta$ satisfies the sixth moment condition~\Ref{thm2.1r}. If $\eta$ is exponentially stabilising in Definition~\ref{defi4r}, then there exist positive constants $\alpha_0$ and $C$ such that 
		\begin{equation}\left|\var\(\bar{W}_{\alpha}\)-\var\(\bar{W}_{\alpha,r}\)\right|\le \frac{1}{\alpha}\label{lma12s3}\end{equation} for all $\alpha\ge\alpha_0$ and $r\ge C\ln (\alpha)$.
	\end{description}
\end{lma}

\noindent{\it Proof.} We start with \Ref{lma12s3}. From Lemma~\ref{lma11}~(b), taking $n=3$ and $k=1$ or $k=3$, we have 
\begin{equation}\max\left\{\|\bar{W}_{\alpha}\|_1,\|\bar{W}_{\alpha,r}\|_1\right\}\le C_0\alpha,~\max\left\{\|\bar{W}_{\alpha}\|_3,\|\bar{W}_{\alpha,r}\|_3\right\}\le C_0\alpha^{\frac{5}{3}}\label{lmapr1}
\end{equation} for some positive constant $C_0\ge 1$. Without loss, we assume $\alpha>1$.
Since \begin{equation}\left|\var\(\bar{W}_{\alpha}\)-\var\(\bar{W}_{\alpha,r}\)\right|\le \left|\mathbb{E}\(\bar{W}_\alpha^2-\bar{W}_{\alpha,r}^2\)\right|+\left|\(\mathbb{E}\bar{W}_\alpha\)^2-\(\mathbb{E}\bar{W}_{\alpha,r}\)^2\right|, \label{lma12.0}
\end{equation} 
using the assumption of stabilisation, we show that each of the terms at the right hand side of \Ref{lma12.0} is bounded by $\frac{1}{2\alpha}$ for $\alpha$ and $r$ sufficiently large. Clearly, the definition of $\bar{W}_{\alpha,r}$ implies that $\bar{W}_\alpha^2-\bar{W}_{\alpha,r}^2= 0$ if $\bar{R}(x,\alpha)\le r$ for all $x\in\overline{\Xi}_{\Gamma_\alpha}$, hence it remains to tackle $E_{r,\alpha}:=\{\bar{R}(x,\alpha)\le r~\mbox{for all }x\in\overline{\Xi}_{\Gamma_\alpha}\}^c$. As shown in the proof of \Ref{lma105.1}, $\mathbb{P}\(E_{r,\alpha}\)\le\alpha C_1e^{-C_2r}$, which, together with \Ref{lmapr1}, H\"older's inequality and Minkowski's inequality,
ensures 
\begin{align}
	\left|\mathbb{E}\(\bar{W}_\alpha^2-\bar{W}_{\alpha,r}^2\)\right|&=\left|\mathbb{E}\[\(\bar{W}_\alpha^2-\bar{W}_{\alpha,r}^2\)\mathbf{1}_{E_{r,\alpha}}\]\right|\nonumber
	\\&\le \|\bar{W}_\alpha^2-\bar{W}_{\alpha,r}^2\|_{\frac{3}{2}}\|\mathbf{1}_{E_{r,\alpha}}\|_{3}\nonumber
	\\&\le \(\|\bar{W}_\alpha^2\|_{\frac{3}{2}}+\|\bar{W}_{\alpha,r}^2\|_{\frac{3}{2}}\)\mathbb{P}(E_{r,\alpha})^{\frac{1}{3}}\nonumber
	\\&= \(\|\bar{W}_\alpha\|_{3}^2+\|\bar{W}_{\alpha,r}\|_{3}^2\)\mathbb{P}(E_{r,\alpha})^{\frac{1}{3}}\le 2C_0^2\alpha^{\frac{10}{3}} \(\alpha C_1e^{-C_2r}\)^{\frac{1}{3}}.\label{lma12.1}
\end{align} 

For the remaining term of \Ref{lma12.0}, we have $$\left|\(\mathbb{E}\bar{W}_\alpha\)^2-\(\mathbb{E}\bar{W}_{\alpha,r}\)^2\right|=\left|\mathbb{E}\bar{W}_\alpha-\mathbb{E}\bar{W}_{\alpha,r}\right|\left|\mathbb{E}\bar{W}_\alpha+\mathbb{E}\bar{W}_{\alpha,r}\right|.$$ The bound \Ref{lmapr1} implies $\left|\mathbb{E}\bar{W}_\alpha+\mathbb{E}\bar{W}_{\alpha,r}\right|\le 2C_0\alpha $. However, using H\"older's inequality, Minkowski's inequality and \Ref{lmapr1} again, we have 
\begin{align}
	\left|\mathbb{E}\bar{W}_\alpha-\mathbb{E}\bar{W}_{\alpha,r}\right|&=\left|\mathbb{E}\[\(\bar{W}_\alpha-\bar{W}_{\alpha,r}\)\mathbf{1}_{E_{r,\alpha}}\]\right|\nonumber
	\\&\le \|\bar{W}_\alpha-\bar{W}_{\alpha,r}\|_3\|\mathbf{1}_{E_{r,\alpha}}\|_{\frac{3}{2}}\nonumber
	\\&\le \(\|\bar{W}_\alpha\|_{3}+\|\bar{W}_{\alpha,r}\|_{3}\)\mathbb{P}(E_{r,\alpha})^{\frac{2}{3}}\nonumber
	\\&\le 2C_0\alpha^{\frac{5}{3}} \(\alpha C_1e^{-C_2r}\)^{\frac{2}{3}},\label{lma12.2}
\end{align}
giving \begin{equation}
	\left|\(\mathbb{E}\bar{W}_\alpha\)^2-\(\mathbb{E}\bar{W}_{\alpha,r}\)^2\right|\le 4C_0^2\alpha^{\frac{8}{3}}\(\alpha C_1e^{-C_2r}\)^{\frac{2}{3}}\label{lma12.20}.
\end{equation} 
We set $r=C\ln(\alpha)$ in the upper bounds of \Ref{lma12.1} and \Ref{lma12.20} and find $C$ such that both bounds are bounded by $1/(2\alpha)$, completing 
the proof of \Ref{lma12s3}.

A line-by-line repetition of the above proof with $\bar{W}_\alpha$ and $\bar{W}_{\alpha,r}$ replaced by $W_\alpha$ and $W_{\alpha,r}$ and $\bar{R}(x,\alpha)$ replaced by $R(x)$ gives \Ref{lma12s1}. \qed

We can now establish the lower bounds for $\var\(W_{\alpha}\)$ and $\var\(\bar{W}_{\alpha}\)$ using the variation conditions \Ref{non-sin} and \Ref{non-sinr}. To this end, we start with a lemma. Recall that $P:=\lambda\mathbb{E}(\eta((\bm{0},M_{\bm{0}}), \Xi_{\bm{0}}+\delta_{(0M_{\bm{0}})})$ for the unrestricted case and $\bar{P}:=\lambda\mathbb{E}(\bar{\eta}((\bm{0},M_{\bm{0}}), \Xi_{\bm{0}}+\delta_{(0,M_{\bm{0}})})$ for the restricted case.

\begin{lma}\label{lma2}
	\begin{description}
		\item{(a)} (unrestricted case) Assume that $\Xi$ satisfies the EDD and the fifth moment condition~\Ref{moment0}, and the score function $\eta$ is translation invariant and satisfies the sixth moment condition~\Ref{thm2.1}. If $\eta$ is exponentially stabilising in Definition~\ref{defi4}, then $$\int_{\Gamma_\alpha}\int_{\mathbb{R}^d}\mathbb{E}\[(g(x,\Xi)\overline{\Xi}(dx)-Pdx)(g(y,\Xi)\overline{\Xi}(dy)-Pdy)\]=\alpha\sigma^2< \infty.$$ Furthermore, for any fixed $\alpha_1>0$, $$\int_{\Gamma_{\alpha_1}}\int_{\mathbb{R}^d\backslash\Gamma_{\alpha_2}}\mathbb{E}\[(g(y,\Xi)\overline{\Xi}(dy)-Pdy)(g(x,\Xi)\overline{\Xi}(dx)-Pdx)\]$$ converges to $0$ exponentially fast as $\alpha_2\rightarrow \infty$.
		\item{(b)} (restricted case)  Assume that $\Xi$ satisfies the EDD and the fifth moment condition~\Ref{moment0}, and the score function $\eta$ is translation invariant and satisfies the sixth moment condition~\Ref{thm2.1r}. If $\eta$ is exponentially stabilising in Definition~\ref{defi4r}, then $$\int_{\Gamma_\alpha}\int_{\mathbb{R}^d}\mathbb{E}\[(\bar{g}(x,\Xi)\overline{\Xi}(dx)-\bar{P}dx)(\bar{g}(y,\Xi)\overline{\Xi}(dy)-\bar{P}dy)\]=\alpha\bar{\sigma}^2< \infty.$$ Furthermore, for any fixed $\alpha_1>0$, $$\int_{\Gamma_{\alpha_1}}\int_{\mathbb{R}^d\backslash\Gamma_{\alpha_2}}\mathbb{E}\[(\bar{g}(y,\Xi)\overline{\Xi}(dy)-\bar{P}dy)(\bar{g}(x,\Xi)\overline{\Xi}(dx)-\bar{P}dx)\]$$ converges to $0$ exponentially fast as $\alpha_2\rightarrow \infty$.
		\end{description}
\end{lma}
\noindent{\it Proof.} \setcounter{con}{3}
We start with the restricted case first. It is sufficient to show that, for any fixed $\alpha_1>0$, $$\int_{\Gamma_{\alpha_1}}\int_{\mathbb{R}^d\backslash\Gamma_{\alpha_2}}\mathbb{E}\[(\bar{g}(y,\Xi)\overline{\Xi}(dy)-\bar{P}dy)(\bar{g}(x,\Xi)\overline{\Xi}(dx)-\bar{P}dx)\]$$ converges to $0$ exponentially fast as $\alpha_2\rightarrow \infty$.   Bearing in mind Remark~\ref{moment1}, at the cost of no more than $C_0(\alpha_1\vee 1)^3$, without loss, we may assume that $\alpha^{1/d}_2>\left(6\alpha^{1/d}_1\right)\vee 2$ and $\frac{\alpha_2^{1/d}}{\ln\left(\alpha_2^{1/d}(2\sqrt{d}+1/3)\right)}\ge 12\theta_3$ with $\theta_3$ in the definition of the EDD. The space $\mathbb{R}^d\backslash \Gamma_{\alpha_2}$ can be divided into sets of the form $\left\{\Gamma_{\alpha_2(1+l)^d}\backslash \Gamma_{\alpha_2l^d}\right\}_{l\in \mathbb{N}}=:\{A_l\}_{l\in \mathbb{N}}$.
Then 
\begin{align*}
	&\int_{\Gamma_{\alpha_1}}\int_{\mathbb{R}^d\backslash\Gamma_{\alpha_2}}\mathbb{E}\[(\bar{g}(y,\Xi)\overline{\Xi}(dy)-\bar{P}dy)(\bar{g}(x,\Xi)\overline{\Xi}(dx)-\bar{P}dx)\]
	\\=&\sum_{l\in \mathbb{N}} \int_{\Gamma_{\alpha_1}}\int_{A_l}\mathbb{E}\[(\bar{g}(y,\Xi)\overline{\Xi}(dy)-\bar{P}dy)(\bar{g}(x,\Xi)\overline{\Xi}(dx)-\bar{P}dx)\]
	\\=&\sum_{l\in \mathbb{N}}\mathbb{E}\[\int_{\Gamma_{\alpha_1}}\int_{A_l}(\bar{g}(y,\Xi)\overline{\Xi}(dy)-\bar{P}dy)(\bar{g}(x,\Xi)\overline{\Xi}(dx)-\bar{P}dx)\]
\end{align*} if the sum in the last line is absolutely convergent.

For $l\in\mathbb{N}$,  $\diam\(B\(A_l,l\alpha_2^{1/d}/6\)\)=\alpha_2^{1/d}(\sqrt{d}(1+l)+l/3)\ge \diam\(B\(\Gamma_{\alpha_1},l\alpha_2^{1/d}/6\)\)=\sqrt{d}\alpha_1^{1/d}+l\alpha_2^{1/d}/3$, and $d\(B\(\Gamma_{\alpha_1},l\alpha_2^{1/d}/6\),B\(A_l,l\alpha_2^{1/d}/6\)\)=\frac{l\alpha_2^{1/d}}{6}-\frac{\alpha_1^{1/d}}{2}\ge \frac{l\alpha_2^{1/d}}{12}$. Since $x/\ln(ax)$ for $a>0$ is an increasing function of $x\ge e/a$, we have
$$\frac{l\alpha_2^{1/d}}{\ln\left(\alpha_2^{1/d}(\sqrt{d}(1+l)+l/3)\right)}\ge\frac{l\alpha_2^{1/d}}{\ln\left(l\alpha_2^{1/d}(2\sqrt{d}+1/3)\right)}\ge \frac{\alpha_2^{1/d}}{\ln\left(\alpha_2^{1/d}(2\sqrt{d}+1/3)\right)}\ge 12\theta_3,$$
which ensures
\begin{align*}
&d\(B\(\Gamma_{\alpha_1},l\alpha_2^{1/d}/6\),B\(A_l,l\alpha_2^{1/d}/6\)\)\\
&\ge \theta_3\ln \(\diam\(B\(\Gamma_{\alpha_1},l\alpha_2^{1/d}/6\)\)\vee \diam\(B\(A_l,l\alpha_2^{1/d}/6\)\)\vee 1\).
\end{align*}
With \Ref{palm4} we can show that $\mathbb{E}\int_A (\bar{g}(x,\Xi)\overline{\Xi}(dx)-\bar{P}dx)=0$ for all bounded measurable set $A\subset \mathbb{R}^d$. Let $\tilde{\Xi}$ denote an independent copy of $\Xi$, then $\mathbb{E}(\int_{\Gamma_{\alpha_1}}(\bar{g}(x,\Xi)\overline{\Xi}(dx)-\bar{P}dx)\int_{A_l}(\bar{g}(y,\tilde{\Xi})\overline{\tilde{\Xi}}(dy)-\bar{P}dy))=0$. For simplicity, we write $S_A:=\int_{A}\bar{g}(x,\Xi)\overline{\Xi}(dx)$, $S_{A,r}:=\int_{A}\bar{g}(x,\Xi)\mathbf{1}_{\bar{R}(x)\le r}\overline{\Xi}(dx)$ and the corresponding counterparts with $\tilde{\Xi}$ instead of $\Xi$ as $\tilde{S}_A$ and $\tilde{S}_{A,r}$ for all bounded measurable sets $A\subset \mathbb{R}^d$.  Using the stabilising condition in Definition~\ref{defi4r} and the EDD, we can get an upper bound for $d_{TV}((S_{\Gamma_{\alpha_1}}, \tilde{S}_{A_l}),(S_{\Gamma_{\alpha_1}}, S_{A_l}))$ as follows:
\begin{align*}
	&d_{TV}((S_{\Gamma_{\alpha_1}}-\alpha_1\bar{P}, \tilde{S}_{A_l}-\Vol(A_l)\bar{P}),(S_{\Gamma_{\alpha_1}}-\alpha_1\bar{P}, S_{A_l}-\Vol(A_l)\bar{P}))
	\\=&d_{TV}((S_{\Gamma_{\alpha_1}}, \tilde{S}_{A_l}),(S_{\Gamma_{\alpha_1}}, S_{A_l}))
	\\\le & \mathbb{P}(S_{\Gamma_{\alpha_1}}\neq S_{\Gamma_{\alpha_1},l\alpha_2^{1/d}/6})+ 2\mathbb{P}(\tilde{S}_{A_l}\neq \tilde{S}_{A_l,l\alpha_2^{1/d}/6})\\
	&+d_{TV}((S_{\Gamma_{\alpha_1},l\alpha_2^{1/d}/6}, \tilde{S}_{A_l,l\alpha_2^{1/d}/6}),(S_{\Gamma_{\alpha_1},l\alpha_2^{1/d}/6}, S_{A_l,l\alpha_2^{1/d}/6}))
	\\\le &\lambda(\alpha_1+2\Vol(A_l))\bar{\tau}(l\alpha_2^{1/d}/6)+\beta_{B(\Gamma_{\alpha_1},l\alpha_2^{1/d}/6),B(A_l,l\alpha_2^{1/d}/6)}
	\\\le& C_{\thecon}\setcounter{xiaa}{\value{con}}\qcon{}e^{-C_{\thecon}l\alpha_2^{1/d}}, \setcounter{xiab}{\value{con}}\qcon{}
\end{align*} for some positive constants $C_{\thexiaa}$ and $C_{\thexiab}$ which are independent of $l$ and $\alpha_2$, where the second inequality follows from the same argument as that for \Ref{lma105.1}, and the last inequality follows from the stabilising condition in Definition~\ref{defi4r} and the EDD. 

Using \cite[p.~254]{BHJ}, we can find a suitable coupling $((X_1,X_2), (Y_1,Y_2))$ of $\(S_{\Gamma_{\alpha_1}}-\alpha_1\bar{P}, \tilde{S}_{A_l}-\Vol(A_l)\bar{P}\)$ and $\(S_{\Gamma_{\alpha_1}}-\alpha_1\bar{P}, S_{A_l}-\Vol(A_l)\bar{P}\)$ such that  $(X_1,X_2)\overset{d}{=}\(S_{\Gamma_{\alpha_1}}-\alpha_1\bar{P}, \tilde{S}_{A_l}-\Vol(A_l)\bar{P}\)$, $(Y_1,Y_2)\overset{d}{=}\(S_{\Gamma_{\alpha_1}}-\alpha_1\bar{P}, S_{A_l}-\Vol(A_l)\bar{P}\)$, $\mathbb{P}(E):=\mathbb{P}((X_1,X_2)\neq(Y_1,Y_2))\le C_{\thexiaa}e^{-C_{\thexiab}l\alpha_2^{1/d}}$. With Remark~\ref{moment2}, we can see that $\|S_{A}-\Vol(A)\bar{P}\|_3\le C_{\thecon}\setcounter{xiae}{\value{con}}\Vol(A)^{\frac{5}{3}}$\qcon{} for $A=\Gamma_{\alpha}$ or $A_l$ for $l\ge 1$, $\alpha>0$ for some constant $C_{\thexiae}$.

From H\"older's inequality and Remark~\ref{moment2}, we can see that 
\begin{align}
	&\left|\mathbb{E}\(\int_{\Gamma_{\alpha_1}}(\bar{g}(x,\Xi)\overline{\Xi}(dx)-\bar{P}dx)\int_{A_l}(\bar{g}(y,\Xi)\overline{\Xi}(dy)-\bar{P}dy)\)\right|
	\nonumber\\=&\left|\mathbb{E}\(\int_{\Gamma_{\alpha_1}}(\bar{g}(x,\Xi)\overline{\Xi}(dx)-\bar{P}dx)\int_{A_l}(\bar{g}(y,\Xi)\overline{\Xi}(dy)-\bar{P}dy)\)\right.
	\nonumber\\&\left.-\mathbb{E}\(\int_{\Gamma_{\alpha_1}}(\bar{g}(x,\Xi)\overline{\Xi}(dx)-\bar{P}dx)\int_{A_l}(\bar{g}(y,\tilde{\Xi})\overline{\tilde{\Xi}}(dy)-\bar{P}dy)\)\right|
	\nonumber\\=&\left|\mathbb{E}\(\int_{\Gamma_{\alpha_1}}(\bar{g}(x,\Xi)\overline{\Xi}(dx)-\bar{P}dx)\int_{A_l}(\bar{g}(y,\Xi)\overline{\Xi}(dy)-\bar{P}dy)\mathbf{1}_{E}\)\right.
	\nonumber\\&\left.-\mathbb{E}\(\int_{\Gamma_{\alpha_1}}(\bar{g}(x,\Xi)\overline{\Xi}(dx)-\bar{P}dx)\int_{A_l}(\bar{g}(y,\tilde{\Xi})\overline{\tilde{\Xi}}(dy)-\bar{P}dy)\mathbf{1}_{E}\)\right|
	\nonumber\\ \le & C_{\thecon}\setcounter{xiac}{\value{con}}\qcon{}e^{-C_{\thecon}\setcounter{xiad}{\value{con}}\qcon{}l\alpha_2^{1/d}}\label{lma2.9-1}
\end{align}  for some positive constants $C_{\thexiac}$ and $C_{\thexiad}$. This, together with the translation invariant property and the definition of $\bar{\sigma}^2$, completes the proof for the restricted case. 

The statement for the unrestricted case can be proved by replacing corresponding counterparts $g$ by $\bar{g}$; $\bar{R}$ by $R$; $\bar{P}$ by $P$; $\bar{\tau}$ by $\tau$ and $\bar{\sigma}^2$ by $\sigma^2$. \qed

To show the volume of the variance, we need a lemma saying that we can approximate $\alpha\sigma^2$ (resp. $\alpha\bar{\sigma}^2$) using the score function $\eta$ restricted to $R\le r$ (resp. $\bar{R}\le r$). For convenience, let $P_r:=\lambda\mathbb{E}(\eta((\bm{0},M_{\bm{0}}), \Xi_{\bm{0}}+\delta_{(\bm{0},M_{\bm{0}})})\mathbf{1}_{R((\bm{0},M_{\bm{0}}), \Xi_{\bm{0}}+\delta_{(\bm{0},M_{\bm{0}})})\le r})$, $\bar{P}_{\alpha,x,r}:=\lambda\mathbb{E}(\eta((x,M_x), \Xi_x+\delta_{(x,M_x)},\Gamma_\alpha)\mathbf{1}_{\bar{R}((x,M_x),\alpha, \Xi_{x}+\delta_{(x,M_x)})\le r})$.

\begin{lma}\label{newlma4}
	\begin{description}
		\item{(a)} (unrestricted case) If the conditions in Lemma~\ref{lma2}~(a) hold and $\sigma^2>0$, then for a fixed sufficiently large $\alpha_0>0$, 
		there exist positive constants $\alpha_1$ and $C$ such that
		\begin{align*}&\mathbb{E}\[\int_{z+\Gamma_{\alpha_0}}(g(x,\Xi)\mathbf{1}_{R(x)\le r}\overline{\Xi}(dx)-P_{r}dx)\int_{\Gamma_\alpha}(g(y,\Xi)\mathbf{1}_{R(y)\le r}\overline{\Xi}(dy)-P_{r}dy)\]\\
		\in& [\frac{1}{2}\alpha_0\sigma^2,\frac{3}{2}\alpha_0\sigma^2]\end{align*} for all $\alpha\ge \alpha_1$ and $r\ge C\ln(\alpha)$, $z\in \Gamma_\alpha$ such that $d(z,\partial \Gamma_\alpha)\ge 5r$.
		\item{(b)} (restricted case)  If the conditions in Lemma~\ref{lma2}~(b) hold and $\bar{\sigma}^2>0$, then for a fixed sufficiently large $\alpha_0>0$, 
		there exist positive constants $\alpha_1$ and $C$ such that
		\begin{align*}&\mathbb{E}\[\int_{z+\Gamma_{\alpha_0}}(g_\alpha(x,\Xi)\mathbf{1}_{\bar{R}(x)\le r}\overline{\Xi}(dx)-\bar{P}_{\alpha,x,r}dx)\int_{\Gamma_\alpha}(g_\alpha(y,\Xi)\mathbf{1}_{\bar{R}(y)\le r}\overline{\Xi}(dy)-\bar{P}_{\alpha,y,r}dy)\]\\
		\in& [\frac{1}{2}\alpha_0\bar{\sigma}^2,\frac{3}{2}\alpha_0\bar{\sigma}^2]\end{align*} for all $\alpha\ge \alpha_1$ and $r\ge C\ln(\alpha)$, $z\in \Gamma_\alpha$ such that $d(z,\partial \Gamma_\alpha)\ge 5r$.
	\end{description}  
\end{lma}
\noindent{\it Proof.} \setcounter{con}{1} 
We prove the restricted case only, and the unrestricted case can be proved similarly. 

Let $\bar{P}_{r}:=\lambda\mathbb{E}(\bar{g}(\bm{0},\Xi_{\bm{0}}+\delta_{(\bm{0},M_{\bm{0}})})\mathbf{1}_{\bar{R}((\bm{0},M_{\bm{0}}),\alpha, \Xi_{\bm{0}}+\delta_{(\bm{0},M_{\bm{0}})})\le r})$, then from the moment condition \Ref{thm2.1r}, we can see that $\max_{\alpha\in \mathbb{R}_+,r\in \mathbb{R}_+,x\in \Gamma_{\alpha}}\{|\bar{P}_{r}|,|\bar{P}_{\alpha,x,r}|\}\le C_{\thecon}$$\setcounter{xiaia}{\thecon}$\qcon{} for some constant $C_{\thexiaia}$. By the translation-invariance, if $x$ and $r$ satisfy $B(x,r)\subset \Gamma_\alpha$, then $g_\alpha(x,\Xi)\mathbf{1}_{\bar{R}(x)\le r}=\bar{g}(x,\Xi)\mathbf{1}_{R(x)\le r}$ and $\bar{P}_r=\bar{P}_{\alpha,x,r}$, where $R(x)=\lim_{\alpha\to\infty}\bar{R}(x)$, see the discussion after Definition~\ref{traninvres0}. The stabilising condition ensures that $\mathbb{P}(\bar{R}((\bm{0},M_{\bm{0}}),\alpha, \Xi_{\bm{0}}+\delta_{(\bm{0},M_{\bm{0}})})\ge r)\le \bar{\tau}(r)$ decreases exponentially fast. Arguing in the same way as that for \Ref{lma12.1}, it follows from the moment condition \Ref{thm2.1r} and H\"older's inequality that $|\bar{P}-\bar{P}_{r}|\le \alpha^{-2}$ for all $r\ge C_{\thecon}\ln(\alpha)$$\setcounter{xiaib}{\thecon}$\qcon{} for some positive constant $C_{\thexiaib}$. 
Writing for brevity $T_{\alpha,r}(x,\Xi,dx):=g_\alpha(x,\Xi)\mathbf{1}_{\bar{R}(x)\le r}\overline{\Xi}(dx)-\bar{P}_{\alpha,x,r}dx$
and $T_r(x,\Xi,dx):=\bar{g}(x,\Xi)\mathbf{1}_{R(x)\le r}\overline{\Xi}(dx)-\bar{P}dx$, we have from Lemma~\ref{lma11}~(b) that there exists a positive constant $\alpha_1'\ge e$ such that 
\begin{align}
&\left|\mathbb{E}\[\int_{z+\Gamma_{\alpha_0}}T_{\alpha,r}(x,\Xi,dx)\int_{B(z+\Gamma_{\alpha_0},3r)}T_{\alpha,r}(y,\Xi,dy)\]\right.\nonumber\\
&\ \ \ -\left.\mathbb{E}\[\int_{z+\Gamma_{\alpha_0}}T_r(x,\Xi,dx)\int_{{B(z+\Gamma_{\alpha_0}},3r)}T_r(y,\Xi,dy)\]\right|\nonumber\\
&\le\frac{1}{8}\alpha_0\bar{\sigma}^2\label{lma5.10-4}\end{align} 
for all $\alpha\ge \alpha_1', r\ge C_{\thexiaib}\ln(\alpha)$.

Following the proof of \Ref{lma105.1}, we have $\mathbb{P}(E^c):=\mathbb{P}\(\{\bar{R}(x)\le r~\mbox{for all}~x\in \overline{\Xi}_{\Gamma_\alpha}\}^c\)\le \alpha\bar{\tau}(r)$. Also, we can see that on the event $E$, $$\int_{z+\Gamma_{\alpha_0}}T_r(x,\Xi,dx)\int_{{B(z+\Gamma_{\alpha_0}},3r)}T_r(y,\Xi,dy)$$ is the same as $$\int_{z+\Gamma_{\alpha_0}}T_\infty(x,\Xi,dx)\int_{{B(z+\Gamma_{\alpha_0}},3r)}T_\infty(y,\Xi,dy),$$ 
where $T_\infty(x,\Xi,dx):=\bar{g}(x,\Xi)\overline{\Xi}(dx)-\bar{P}dx$. Using H\"older's inequality as for \Ref{lma12.1}, we can find a $C_{\thecon}\in\mathbb{R}_+$\setcounter{xiaic}{\thecon}\qcon{}  such that
\begin{align}
&\left|\int_{z+\Gamma_{\alpha_0}}T_r(x,\Xi,dx)\int_{{B(z+\Gamma_{\alpha_0}},3r)}T_r(y,\Xi,dy)\right.\nonumber\\
&\ \ \ \left.-\int_{z+\Gamma_{\alpha_0}}T_\infty(x,\Xi,dx)\int_{{B(z+\Gamma_{\alpha_0}},3r)}T_\infty(y,\Xi,dy)\right|\nonumber\\
&\le \frac{1}{8}\alpha_0\bar{\sigma}^2\label{lma5.10-5}
\end{align}
for all $r\ge C_{\thexiaic}\ln(\alpha)$.

Using the translation invariant property, and replacing $\mathbb{R}^d\backslash \Gamma_{\alpha_2}$ by $\mathbb{R}^d\backslash B(\Gamma_{\alpha_0},3r)$ in the proof of Lemma~\ref{lma2}, with necessary minor adjustment,  we obtain
\begin{align}
	&\left|\mathbb{E}\[\int_{z+\Gamma_{\alpha_0}}T_\infty(x,\Xi,dx)\int_{{B(z+\Gamma_{\alpha_0}},3r)}T_\infty(y,\Xi,dy)\]-\alpha_0\bar{\sigma}^2\right|\nonumber
	\\ &=\left|\mathbb{E}\[\int_{\Gamma_{\alpha_0}}T_\infty(x,\Xi,dx)\int_{B(\Gamma_{\alpha_0},3r)}T_\infty(y,\Xi,dy)\]-\alpha_0\bar{\sigma}^2\right|\nonumber\\
	&\le \frac{1}{8}\alpha_0\bar{\sigma}^2\label{lma5.10-1}
\end{align}
for $r\ge C_{\thecon}\ln(\alpha)$\setcounter{xiaie}{\thecon}\qcon{}. 
On the other hand, the EDD ensures that we can find an independent copy $\tilde{\Xi}$ of $\Xi$ such that 
$$E_1:=\{\Xi_{B(z+\Gamma_{\alpha_0},r)\cup (\Gamma_{\alpha}\backslash B(z+\Gamma_{\alpha_0},2r)) }\ne \Xi_{B(z+\Gamma_{\alpha_0},r)}\cup \tilde{\Xi}_{\Gamma_{\alpha}\backslash B(z+\Gamma_{\alpha_0},2r) }\}$$
satisfies 
\begin{equation}\mathbb{P}(E_1)=\beta_{B(z+\Gamma_{\alpha_0},r), \Gamma_{\alpha}\backslash B(z+\Gamma_{\alpha_0},2r)}\le \theta_1\(\(\sqrt{d}\alpha^{2\theta_0/d}\)\vee 1\)
e^{-\theta_2r}.\label{lma5.10-2}\end{equation} Since $\mathbb{E}T_{\alpha,r}(x,\Xi,dx)=0$, following the same argument as that for \Ref{lma2.9-1} and applying H\"older's inequality in the first inequality, the moment condition \Ref{thm2.1r} and \Ref{lma5.10-2} in the last inequality below, we get
\begin{align}
&\left|\mathbb{E}\[\int_{z+\Gamma_{\alpha_0}}T_{\alpha,r}(x,\Xi,dx)\int_{\Gamma_{\alpha}\backslash B(z+\Gamma_{\alpha_0},3r)}T_{\alpha,r}(y,\Xi,dy)\]\right|\nonumber\\
&=\left|\mathbb{E}\[{\bf 1}_{E_1}\int_{z+\Gamma_{\alpha_0}}T_{\alpha,r}(x,\Xi,dx)\int_{\Gamma_{\alpha}\backslash B(z+\Gamma_{\alpha_0},3r)}T_{\alpha,r}(y,\Xi,dy)\]\right.\nonumber\\
&\ \ \ -{\bf 1}_{E_1}\left.\mathbb{E}\[\int_{z+\Gamma_{\alpha_0}}T_{\alpha,r}(x,\Xi,dx)\int_{\Gamma_{\alpha}\backslash B(z+\Gamma_{\alpha_0},3r)}T_{\alpha,r}(y,\tilde{\Xi},dy)\]\right|\nonumber\\
&\le \left\|\int_{z+\Gamma_{\alpha_0}}T_{\alpha,r}(x,\Xi,dx)\int_{\Gamma_{\alpha}\backslash B(z+\Gamma_{\alpha_0},3r)}T_{\alpha,r}(y,\Xi,dy)\right\|_{3/2}\mathbb{P}(E_1)^{1/3}\nonumber\\
&\ \ \ +\left\|\int_{z+\Gamma_{\alpha_0}}T_{\alpha,r}(x,\Xi,dx)\int_{\Gamma_{\alpha}\backslash B(z+\Gamma_{\alpha_0},3r)}T_{\alpha,r}(y,\tilde{\Xi},dy)\right\|_{3/2}\mathbb{P}(E_1)^{1/3}\nonumber\\
&\le \frac{1}{8}\alpha_0\bar{\sigma}^2\label{lma5.10-3}
\end{align}
for all $z$ satisfying $d(z,\partial \Gamma_\alpha)\ge 5r$, $r\ge C_{\thecon}\ln(\alpha)$\setcounter{xiaid}{\thecon} and $\alpha\ge \alpha_2'$.  Collecting \Ref{lma5.10-4}, \Ref{lma5.10-5}, \Ref{lma5.10-1} and \Ref{lma5.10-3}, we obtain claim (b) with $C=\max\{C_{\thexiaib}, C_{\thexiaic},C_{\thexiaie},C_{\thexiaid}\}$ and $\alpha_1=\max\{\alpha_1',\alpha_2'\}$.

The statement for the unrestricted case can be proved by replacing $\bar{P}_{\alpha,x,r}$ and $\bar{P}_r$ by $P_{\alpha,x,r}$ and $P_r$; $\bar{g}$ by $g$; $\bar{R}$ by $R$; $\bar{\tau}$ by $\tau$; $\bar{P}$ by $P$.
\qed

Together with the variation conditions \Ref{non-sin} and \Ref{non-sinr}, we can show that the variances of $W_\alpha$ and $\bar{W}_{\alpha}$ are of the order $\alpha$ as claimed in Theorem~\ref{thmvar}.

\noindent{\it Proof of Theorem~\ref{thmvar}.} We show the statement for the restricted case first, and the statement for the unrestricted case can be shown in the same way. 

\setcounter{con}{1}
To begin with, we choose $\alpha_0$ and $C$ such that Lemma~\ref{lma12} holds, i.e., $|\var(\bar{W}_\alpha)-\var(\bar{W}_{\alpha,r})|\le \frac1\alpha$ for $r\ge C\ln(\alpha)$ and $\alpha\ge \alpha_0$. We then choose $C_{\thecon}\ge C$\setcounter{xiaha}{\thecon}\qcon{} and $\alpha_1\ge \alpha_0\vee e $ such that Lemma~\ref{newlma4}~(b) holds for all $r\ge C_{\thexiaha}\ln(\alpha)$ and $\alpha\ge \alpha_1$. In the rest of the proof, we fix $r= C_{\thexiaha}\ln(\alpha)$. 
Replacing $\theta_3$ in the definition of the EDD with $2C_{\thexiaha}$ if necessary, we assume $\theta_3\ge 2C_{\thexiaha}$. Cover $\Gamma_\alpha\backslash B(\partial \Gamma_\alpha, 5\theta_3\ln(\alpha))$ with disjoint cubes $\mathbb{C}_{1},\dots, \mathbb{C}_{n_\alpha}$ each with a volume between $\alpha_0$ and $2\alpha_0$ and intersects $\Gamma_\alpha\backslash B(\partial \Gamma_\alpha, 5\theta_3\ln(\alpha))$, then the number of cubes $n_\alpha$ has the same order as $\alpha$. For convenience, let $\mathbb{C}_\alpha:=\cup_{1\le i\le n_\alpha} \mathbb{C}_i$, then from Lemma~\ref{newlma4}~(b), for $\alpha\ge \alpha_1$, 
\begin{align}
&\mathbb{E}\[\int_{\mathbb{C}_\alpha}(g_\alpha(x,\Xi)\mathbf{1}_{\bar{R}(x)\le r}\overline{\Xi}(dx)-\bar{P}_{\alpha,x,r}dx)\int_{\Gamma_\alpha}(g_\alpha(y,\Xi)\mathbf{1}_{\bar{R}(y)\le r}\overline{\Xi}(dy)-\bar{P}_{\alpha,y,r}dy)\]\nonumber\\
&\in [\frac{1}{2}n_\alpha\alpha_0\bar{\sigma}^2,3n_\alpha\alpha_0\bar{\sigma}^2].\label{newlma3.0}
\end{align}
Next, we show that
\begin{align}&\left|\mathbb{E}\[\int_{\Gamma_\alpha\backslash\mathbb{C}_\alpha}(g_\alpha(x,\Xi)\mathbf{1}_{\bar{R}(x)\le r}\overline{\Xi}(dx)-\bar{P}_{\alpha,x,r}dx)\int_{\Gamma_\alpha}(g_\alpha(y,\Xi)\mathbf{1}_{\bar{R}(y)\le r}\overline{\Xi}(dy)-\bar{P}_{\alpha,y,r}dy)\]\right|\nonumber\\
&\le O\((\ln \alpha)^{2d+1}\alpha^{\frac{d-1}{d}}\).\label{proof2.11xiaa1}
\end{align}
To this end, we choose $C_{\thecon}>1\setcounter{xiahb}{\thecon}$\qcon{} such that $\theta_2\theta_3C_{\thexiahb}\ge 13.5+\theta_0/d$ and $C_{\thexiahb}\theta_3>2C_{\thexiaha}$, where $\theta_0,\theta_2$ are as in the EDD, divide $\Gamma_\alpha\backslash\mathbb{C}_\alpha$ into at most $n_{\mathbb{D}}=O\(\alpha^{\frac{d-1}{d}}\ln (\alpha)\)$ cubes $\{\mathbb{D}_i\}_{1\le i\le n_\mathbb{D} }$ having diameters between $1$ and $\theta_3\ln(\alpha)$. For each $x\in\mathbb{D}_i$, $T(x,\Xi,dx):=g_\alpha(x,\Xi)\mathbf{1}_{\bar{R}(x)\le r}\overline{\Xi}(dx)-\bar{P}_{\alpha,x,r}dx$ is completely determined by $\Xi_{B(\mathbb{D}_i,r)}$ and $T(y,\Xi,dy)$ is completely determined by $\Xi_{\Gamma_\alpha\cap B(\mathbb{D}_i,2\theta_3C_{\thexiahb}\ln(\alpha))^c}$ if $x\in B(\mathbb{D}_i,r)$ and $\|y-x\|>4\theta_3C_{\thexiahb}\ln(\alpha)$. By the EDD, we can find an independent copy $\tilde{\Xi}$ of $\Xi$ such that
\begin{align}
\mathbb{P}(E_i):=&\mathbb{P}\(\Xi_{B(\mathbb{D}_i,r)}\cup\Xi_{\Gamma_\alpha\cap B(\mathbb{D}_i,2\theta_3C_{\thexiahb}\ln(\alpha))^c}\ne
\Xi_{B(\mathbb{D}_i,r)}\cup\tilde{\Xi}_{\Gamma_\alpha\cap B(\mathbb{D}_i,2\theta_3C_{\thexiahb}\ln(\alpha))^c}\)\nonumber\\
=&\beta_{B(\mathbb{D}_i,r),\Gamma_\alpha\cap B(\mathbb{D}_i,2\theta_3C_{\thexiahb}\ln(\alpha))^c}
\le C_{\thecon}\ln(\alpha)^{\theta_0}\alpha^{\theta_0/d}e^{-\theta_2\theta_3C_{\thexiahb}\ln(\alpha)}.\label{proof2.11xiaa2}
\end{align}\qcon{}
Hence, following the same steps as those for \Ref{lma2.9-1} and using H\"older's inequality in the first inequality, \Ref{proof2.11xiaa2} and Remark~\ref{moment1} with $k=3$  in the second inequality below, we have
\begin{align}&\left|\mathbb{E}\[\int_{x\in\mathbb{D}_i}T(x,\Xi,dx)\int_{y\in\Gamma_\alpha,\|y-x\|>4\theta_3C_{\thexiahb}\ln(\alpha)}T(y,\Xi,dy)\]\right|\nonumber\\
\ignore{&=\left|\mathbb{E}\[\int_{x\in\mathbb{D}_i}T(x,\Xi,dx)\int_{y\in\Gamma_\alpha,\|y-x\|>4\theta_3C_{\thexiahb}\ln(\alpha)}T(y,\Xi,dy)\]\right.\nonumber\\
&\ \ \ \ \ \ \  -\left.\mathbb{E}\[\int_{x\in\mathbb{D}_i}T(x,\Xi,dx)\int_{y\in\Gamma_\alpha,\|y-x\|>4\theta_3C_{\thexiahb}\ln(\alpha)}T(y,\tilde{\Xi},dy)\]\right|\nonumber\\}
&=\left|\mathbb{E}\[{\bf 1}_{E_i}\int_{x\in\mathbb{D}_i}T(x,\Xi,dx)\int_{y\in\Gamma_\alpha,\|y-x\|>4\theta_3C_{\thexiahb}\ln(\alpha)}T(y,\Xi,dy)\]\right.\nonumber\\
&\ \ \ \ \ \ \ -\left.\mathbb{E}\[{\bf 1}_{E_i}\int_{x\in\mathbb{D}_i}T(x,\Xi,dx)\int_{y\in\Gamma_\alpha,\|y-x\|>4\theta_3C_{\thexiahb}\ln(\alpha)}T(y,\tilde{\Xi},dy)\]\right|\nonumber\\
&\le \left\|\int_{x\in\mathbb{D}_i}T(x,\Xi,dx)\int_{y\in\Gamma_\alpha,\|y-x\|>4\theta_3C_{\thexiahb}\ln(\alpha)}T(y,\Xi,dy)\right\|_{3/2}\mathbb{P}(E_i)^{1/3}\nonumber\\
&\ \ \ \ \ \ \ +\left\|\int_{x\in\mathbb{D}_i}T(x,\Xi,dx)\int_{y\in\Gamma_\alpha,\|y-x\|>4\theta_3C_{\thexiahb}\ln(\alpha)}T(y,\tilde{\Xi},dy)\right\|_{3/2}\mathbb{P}(E_i)^{1/3}\nonumber\\
&\le O\(\alpha^{\frac52+\frac{\theta_0}{3d}-\frac13\theta_2\theta_3C_{\thexiahb}}\ln(\alpha)^{\theta_0/3}\)\le O\(\ln(\alpha)^{\theta_0/3}\alpha^{-2}\).\label{proof2.11xiaa3}
\end{align}
Adding the estimates of \Ref{proof2.11xiaa3} for $1\le i\le n_{\mathbb{D}}$ and using the fact $n_{\mathbb{D}}=O\(\alpha^{\frac{d-1}{d}}\ln \alpha\)$, we obtain
\begin{align}\left|\mathbb{E}\[\int_{x\in\Gamma_\alpha\backslash\mathbb{C}_\alpha}T(x,\Xi,dx)\int_{y\in\Gamma_\alpha,\|y-x\|>4\theta_3C_{\thexiahb}\ln(\alpha)}T(y,\Xi,dy)\]\right|
\le O\(\alpha^{-1}\).\label{proof2.11xiaa4}
\end{align}
For the remaining part, we have
\begin{align}&\left|\mathbb{E}\[\iint_{x\in\Gamma_\alpha\backslash\mathbb{C}_\alpha,y\in\Gamma_\alpha,\|y-x\|\le 4\theta_3C_{\thexiahb}\ln(\alpha)}g_\alpha(x,\Xi)\mathbf{1}_{\bar{R}(x)\le r}\overline{\Xi}(dx)g_\alpha(y,\Xi)\mathbf{1}_{\bar{R}(y)\le r}\overline{\Xi}(dy)\]\right|\nonumber\\
&\le\mathbb{E}\[\iint_{x\in\Gamma_\alpha\backslash\mathbb{C}_\alpha,y\in\Gamma_\alpha,\|y-x\|\le 4\theta_3C_{\thexiahb}\ln(\alpha)}\frac12\(g_\alpha(x,\Xi)^2\mathbf{1}_{\bar{R}(x)\le r}+g_\alpha(y,\Xi)^2\mathbf{1}_{\bar{R}(y)\le r}\)\overline{\Xi}(dx)\overline{\Xi}(dy)\]\nonumber\\
&\le\frac12\mathbb{E}\[\int_{x\in\Gamma_\alpha\backslash\mathbb{C}_\alpha}\overline{\Xi}(B(x,4\theta_3C_{\thexiahb}\ln(\alpha)))g_\alpha(x,\Xi)^2\mathbf{1}_{\bar{R}(x)\le r}\overline{\Xi}(dx)\]\nonumber\\ 
&\ \ \  +\frac12\mathbb{E}\[\int_{y\in\Gamma_\alpha\cap B(\partial \Gamma_\alpha,9\theta_3C_{\thexiahb}\ln(\alpha))}\overline{\Xi}(B(y,4\theta_3C_{\thexiahb}\ln(\alpha)))g_\alpha(y,\Xi)^2\mathbf{1}_{\bar{R}(y)\le r}\overline{\Xi}(dy)\]\nonumber\\
&\le\mathbb{E}\[\int_{y\in\Gamma_\alpha\cap B(\partial \Gamma_\alpha,9\theta_3C_{\thexiahb}\ln(\alpha))}\overline{\Xi}(B(y,4\theta_3C_{\thexiahb}\ln(\alpha)))g_\alpha(y,\Xi)^2\mathbf{1}_{\bar{R}(y)\le r}\overline{\Xi}(dy)\]\nonumber\\
&\le \mathbb{E}\[\int_{y\in\Gamma_\alpha\cap B(\partial \Gamma_\alpha,9\theta_3C_{\thexiahb}\ln(\alpha))}
\frac12\(\overline{\Xi}(B(y,4\theta_3C_{\thexiahb}\ln(\alpha)))^2+g_\alpha(y,\Xi)^4\mathbf{1}_{\bar{R}(y)\le r}\)\overline{\Xi}(dy)\]\nonumber\\
&\le \int_{y\in\Gamma_\alpha\cap B(\partial \Gamma_\alpha,9\theta_3C_{\thexiahb}\ln(\alpha))}
\(\(\overline{\Xi}_{\bf 0}(B({\bf 0},4\theta_3C_{\thexiahb}\ln(\alpha)))+1\)^2\lambda dy+g_\alpha(y,\Xi)^4\mathbf{1}_{\bar{R}(y)\le r}\overline{\Xi}(dy)\).\label{proof2.11xiaa5}
\end{align}
However, by Remark~\ref{remoment}, 
\begin{align*}
&O\(\ln(\alpha)^{3d}\)\\=&\mathbb{E}\overline{\Xi}(B({\bf 0},8\theta_3C_{\thexiahb}\ln(\alpha)))^3\nonumber\\
\ge &\mathbb{E}\int_{B({\bf 0},4\theta_3C_{\thexiahb}\ln(\alpha))}\overline{\Xi}(B(x,4\theta_3C_{\thexiahb}\ln(\alpha)))^2\overline{\Xi}(dx)\nonumber\\
=&\lambda\int_{B({\bf 0},4\theta_3C_{\thexiahb}\ln(\alpha))}\mathbb{E}\(\overline{\Xi}_x(B(x,4\theta_3C_{\thexiahb}\ln(\alpha)))+1\)^2dx\nonumber\\
=&\frac{\lambda (4\theta_3C_{\thexiahb}\ln(\alpha))^d\pi^{d/2}}{\Gamma(1+d/2)}\mathbb{E}\(\overline{\Xi}_{\bf 0}(B({\bf 0},4\theta_3C_{\thexiahb}\ln(\alpha)))+1\)^2,
\end{align*}
which implies
\begin{equation}
\mathbb{E}\(\overline{\Xi}_{\bf 0}(B({\bf 0},4\theta_3C_{\thexiahb}\ln(\alpha)))+1\)^2\le O\(\ln(\alpha)^{2d}\).\label{proof2.11xiaa6}
\end{equation}
Combining \Ref{thm2.1r}, \Ref{proof2.11xiaa5} and \Ref{proof2.11xiaa6} gives
\begin{align}&\left|\mathbb{E}\[\iint_{x\in\Gamma_\alpha\backslash\mathbb{C}_\alpha,y\in\Gamma_\alpha,\|y-x\|\le 4\theta_3C_{\thexiahb}\ln(\alpha)}g_\alpha(x,\Xi)\mathbf{1}_{\bar{R}(x)\le r}\overline{\Xi}(dx)g_\alpha(y,\Xi)\mathbf{1}_{\bar{R}(y)\le r}\overline{\Xi}(dy)\]\right|\nonumber\\
&\le O\((\ln \alpha)^{2d+1}\alpha^{\frac{d-1}{d}}\). \label{proof2.11xiaa7}
\end{align}
Direct verification using \Ref{thm2.1r} again gives
\begin{align}\left|\mathbb{E}\[\iint_{x\in\Gamma_\alpha\backslash\mathbb{C}_\alpha,y\in\Gamma_\alpha,\|y-x\|\le 4\theta_3C_{\thexiahb}\ln(\alpha)}\bar{P}_{\alpha,x,r}dx\bar{P}_{\alpha,y,r}dy\]\right|
\le O\(\alpha^{\frac{d-1}{d}}\ln(\alpha)^{d+1}\). \label{proof2.11xiaa8}
\end{align}
Collecting \Ref{proof2.11xiaa4}, \Ref{proof2.11xiaa7} and \Ref{proof2.11xiaa8}, we have \Ref{proof2.11xiaa1}. 

Now, since $\bar W_{\alpha,r}$ can be decomposed as
\begin{align*}
	&\var(\bar W_{\alpha,r})\nonumber\\
%	=&\mathbb{E}\[(\bar W_{\alpha,r}-\mathbb{E}(\bar W_{\alpha,r}))^2\]\nonumber\\
%	=&\mathbb{E}\[\int_{\Gamma_\alpha}(g_\alpha(x,\Xi)\mathbf{1}_{\bar{R}(x)\le r}\overline{\Xi}(dx)-\bar{P}_{\alpha,x,r}dx)\int_{\Gamma_\alpha}(g_\alpha(y,\Xi)\mathbf{1}_{\bar{R}(y)\le r}\overline{\Xi}(dy)-\bar{P}_{\alpha,y,r}dy)\]\nonumber
	\\ =&\mathbb{E}\[\int_{\mathbb{C}_\alpha}(g_\alpha(x,\Xi)\mathbf{1}_{\bar{R}(x)\le r}\overline{\Xi}(dx)-\bar{P}_{\alpha,x,r}dx)\int_{\Gamma_\alpha}(g_\alpha(y,\Xi)\mathbf{1}_{\bar{R}(y)\le r}\overline{\Xi}(dy)-\bar{P}_{\alpha,y,r}dy)\]\nonumber
	\\&+\mathbb{E}\[\int_{\Gamma_\alpha\backslash\mathbb{C}_\alpha}(g_\alpha(x,\Xi)\mathbf{1}_{\bar{R}(x)\le r}\overline{\Xi}(dx)-\bar{P}_{\alpha,x,r}dx)\int_{\Gamma_\alpha}(g_\alpha(y,\Xi)\mathbf{1}_{\bar{R}(y)\le r}\overline{\Xi}(dy)-\bar{P}_{\alpha,y,r}dy)\],%\label{newlma3.1}
\end{align*}
it follows from \Ref{newlma3.0} and \Ref{proof2.11xiaa1} that, for $\alpha\ge \alpha_1\vee \alpha_0\vee e$, $\var(\bar{W}_{\alpha,r})\in [C_{\thecon}\setcounter{xiahc}{\thecon}\qcon{}\alpha, C_{\thecon}\setcounter{xiahd}{\thecon}\alpha]$\qcon{} for some positive constants $C_{\thexiahc}$, $C_{\thexiahd}$. This, together with $|\var(\bar{W}_\alpha)-\var(\bar{W}_{\alpha,r})|\le \frac1\alpha$, ensures that $\var(\bar{W}_\alpha)=\Theta(\var(\bar{W}_{\alpha,r}))=\Theta(\alpha)$.

The statement for the unrestricted case can be proved by replacing $g_\alpha$ by $g$; $\bar{P}_{\alpha,x,r}$ by $P_r$; $\bar{W}_\alpha$ by $W_\alpha$ and $\bar{W}_{\alpha,r}$ by $W_{\alpha,r}$. \qed
\begin{re}\label{re1}
	Since the variance is always non-negative, the proof of this also shows that $\sigma^2$ (resp. $\bar{\sigma}^2$) defined in \Ref{non-sin} (resp. \Ref{non-sinr}) is non-negative.
\end{re}

\begin{re}
	Using the same idea as in the proof of Theorem~\ref{thmvar}, we can see that $\var(W_{\alpha})$ and $\var(\bar{W}_\alpha)$ cannot have an order greater than $\alpha$.
\end{re}

\noindent{\it Proof of Theorem~\ref{thm2.a2}.}  \setcounter{con}{1}
Let $\mu_\alpha:= \mathbb{E}\(\bar{W}_\alpha\)$, $\mu_{\alpha,r}:=\mathbb{E}\(\bar{W}_{\alpha,r}\)$, $\sigma_{\alpha}^2:=\var\(\bar{W}_{\alpha}\)$, $\sigma_{\alpha,r}^2:=\var\(\bar{W}_{\alpha,r}\)$ and
$\bar{Z}_{\alpha,r}\sim N\(\frac{\mu_{\alpha,r}-\mu_\alpha}{\sigma_\alpha}, \frac{\sigma_{\alpha,r}^2}{\sigma_\alpha^2}\)$, then it follows from the triangle inequality that
\begin{equation}d_W\left(\frac{\bar{W}_\alpha-\mu_{\alpha}}{\sigma_{\alpha}},Z\right)\le d_W\(\frac{\bar{W}_\alpha-\mu_{\alpha}}{\sigma_{\alpha}},\frac{\bar{W}_{\alpha,r}-\mu_{\alpha}}{\sigma_\alpha}\)+d_W\(Z,\bar{Z}_{\alpha,r}\)+d_W\(\frac{\bar{W}_{\alpha,r}-\mu_{\alpha}}{\sigma_\alpha},\bar{Z}_{\alpha,r}\).\label{thm2a01}
\end{equation} Next, we analyse the three terms of \Ref{thm2a01} separately. We start with the exponentially stabilising case (ii).   

The first term of \Ref{thm2a01} can be bounded using Lemma~\ref{lma105}~(b), the variation condition and the property of the Wasserstein distance. Let $U_{\alpha}:=\(\bar{W}_{\alpha}-\mu_{\alpha}\)/\sigma_{\alpha}$ and $U_{\alpha,r}:=\(\bar{W}_{\alpha,r}-\mu_{\alpha,r}\)/\sigma_{\alpha}$. According to the property of the total variation distance and \cite[p.~254]{BHJ}, we can find a coupling $(\bar{U}_{\alpha},\bar{U}_{\alpha,r})$ of $U_{\alpha}$ and $U_{\alpha,r}$ such that $\bar{U}_{\alpha}\overset{d}{=}U_{\alpha}$,  $\bar{U}_{\alpha,r}\overset{d}{=}U_{\alpha,r}$ and $$\mathbb{P}(\bar{U}_{\alpha}\neq \bar{U}_{\alpha,r})=:\mathbb{P}(E_{\alpha,r})=d_{TV}(U_{\alpha},U_{\alpha,r})=d_{TV}(\bar{W}_{\alpha},\bar{W}_{\alpha,r})\le C_{\thecon}\qcon{}\alpha e^{-C_{\thecon}\qcon{}r}.$$
Then from H\"older's inequality, the variation condition and Lemma~\ref{lma11}, 
	\begin{align}
		&d_W\(\frac{\bar{W}_\alpha-\mu_{\alpha}}{\sigma_{\alpha}},\frac{\bar{W}_{\alpha,r}-\mu_{\alpha}}{\sigma_\alpha}\)\nonumber
		\\=&\inf_{X\overset{d}{=}U_{\alpha},Y\overset{d}{=}U_{\alpha,r}}\mathbb{E}(|X-Y|)\nonumber 
			\\\le& \mathbb{E}(|\bar{U}_{\alpha}-\bar{U}_{\alpha,r}|)\le \mathbb{E}\(|\bar{U}_{\alpha}|\mathbf{1}_{E_{\alpha,r}}\)+\mathbb{E}\(|\bar{U}_{\alpha,r}|\mathbf{1}_{E_{\alpha,r}}\)\le \frac1\alpha,\label{thm2a02}
	\end{align}
for $r>C_{\thecon}\ln(\alpha).$ \setcounter{cproofa}{\thecon} \qcon{}

For the second term of \Ref{thm2a01}, we can establish an upper bound using Lemma~\ref{lma10}. To this end, \Ref{lma12s3} gives
\begin{equation}
	\left|\sigma_\alpha^2-\sigma_{\alpha,r}^2\right|\le \frac1\alpha,
	\label{thm2a04}
\end{equation}
which, together with the condition given in the theorem, implies
\begin{equation}
	\sigma_{\alpha,r}^2= \Omega(\alpha^\nu), \ \ \ \ \ \sigma_{\alpha}^2= \Omega(\alpha^\nu),
	\label{thm2a03}
\end{equation}
for $r>C_{\thecon}\ln(\alpha).$\setcounter{cproofc}{\thecon}\qcon{}
 We combine \Ref{lma12.2} and \Ref{thm2a03}  to obtain
 \begin{equation}\frac{\left|\mu_\alpha-\mu_{\alpha,r}\right|}{\sigma_\alpha}\le O\(\alpha^{-1}\),
 	\label{thm2a14}
 \end{equation}
 for $r>C_{\thecon} \qcon{}\ln(\alpha).$
 Therefore, it follows from \Ref{thm2a04}, \Ref{thm2a03}, \Ref{thm2a14} and Lemma~\ref{lma10} that
 \begin{align}
 	d_{W}(Z,\bar{Z}_{\alpha,r})
 	\le\frac{\left|\mu_\alpha-\mu_{\alpha,r}\right|}{\sigma_\alpha}+\frac{\left|\sigma_{\alpha,r}-\sigma_{\alpha}\right|}{\sigma_{\alpha}}\le O(\alpha^{-1})\label{thm2a05}
 \end{align} for $r>C_{\thecon} \ln(\alpha).$\setcounter{cproofb}{\thecon}

It remains to tackle the last term of \Ref{thm2a01}. From the definition of the Wasserstein distance, we have 
\begin{align}
	&d_W\(\frac{\bar{W}_{\alpha,r}-\mu_{\alpha}}{\sigma_\alpha},\bar{Z}_{\alpha,r}\)= d_W\(\frac{\bar{W}_{\alpha,r}-\mu_{\alpha,r}}{\sigma_\alpha},\bar{Z}_{\alpha,r}+\frac{\mu_\alpha-\mu_{\alpha,r}}{\sigma_\alpha}\) \nonumber
	\\&\le\frac{\sigma_{\alpha,r}}{\sigma_{\alpha}}d_{W}\(\frac{\bar{W}_{\alpha,r}-\mu_{\alpha,r}}{\sigma_{\alpha,r}},Z\) \le 2d_{W}\(V_{\alpha,r},Z\)\label{distance1}
\end{align} for $r>C_{\thecproofc}\ln(\alpha)$ when $\alpha$ large, where $V_{\alpha,r}:=\(\bar{W}_{\alpha,r}-\mu_{\alpha,r}\)/\sigma_{\alpha,r}$.
We now use Stein's method to bound the Wasserstein distance between $V_{\alpha,r}$ and $Z$. Stein's method for the normal approximation hinges on a Stein equation (see  \cite[pp.~15--16]{CGS11})
\begin{equation}f'(w)-wf(w)=h(w)-Nh,\label{steineq1}\end{equation} where $Nh:=\mean h(Z)$. The solution of \Ref{steineq1} is given by $$f_h(w)=e^{w^2/2}\int^w_{-\infty}e^{-t^2/2}(h(t)-Nh)dt=-e^{w^2/2}\int_{w}^{\infty}e^{-t^2/2}(h(t)-Nh)dt.$$ Recall the definition of the Wasserstein distance \Ref{defWass}, we have, for any random variable $X$, 
\begin{equation}\label{Stein}
d_W(X,Z)=\sup_{h\in \mathscr{F}_{\rm Lip}}|\mathbb{E}(h(X)-h(Z))|\le\sup_{f\in \mathscr{F}}\left|\mathbb{E}\(f'(X)-Xf(X)\)\right|,
\end{equation} where $\mathscr{F}:=\left\{f;~\mathbb{R}\rightarrow\mathbb{R}, \|f\|\le 2,~\|f'\|\le \sqrt{\frac{2}{\pi}},~\|f''\|\le 2\right\}$. From the definition of $V_{\alpha,r}$, we can represent it as $V_{\alpha,r}=\frac1{\sigma_{\alpha,r}}\int_{\Gamma_\alpha} \(g_\alpha(x,\Xi){\mathbf1}_{\bar{R}(x)\le r}\overline{\Xi}(dx)-P_{\alpha,x,r}dx\)=:\int_{\Gamma_\alpha}V(dx)$ if this does not cause confusion. Then, from the definition of $V_{\alpha,r}$, we have 
\begin{align*}
	1=\var(V_{\alpha,r})%=\sigma_{\alpha,r}^{-2}\E\(\int_{\Gamma_\alpha}g_\alpha(x,\Xi){\mathbf1}_{\bar{R}(x)\le r}\overline{\Xi}(dx)-\bar{P}_{\alpha,x,r}dx\)^2
	=\E\(\int_{\Gamma_\alpha}V(dx)\)^2. %\label{variance}
\end{align*}
To bound $d_W(V_{\alpha,r},Z)$, by \Ref{Stein}, it is sufficient to bound \begin{equation}\label{Stein1}|\mathbb{E}\(f'(V_{\alpha,r})-V_{\alpha,r}f(V_{\alpha,r})\)|
\end{equation} for all $f\in \mathscr{F}$. To do this, let's consider the two terms separately. 

Before analysing \Ref{Stein1}, for large $\alpha$, we can divide $\Gamma_\alpha$ into disjoint cubes with volumes at least $\alpha_0\le \alpha$ for some positive constant $\alpha_0$. To this end, we can find a partition of $\Gamma_\alpha$, $\mathscr{C}:=\{\mathbb{C}_1,\dots,\mathbb{C}_{n_\alpha}\}$, where $\mathbb{C}_i$ are cubes with edge length $\frac{\alpha^{1/d}}{\left\lfloor(\alpha/\alpha_0)^{1/d}\right\rfloor}$ for all $1\le i\le n_{\alpha}$. Then each cube in $\mathscr{C}$ has a volume no more than $2^d\alpha_0$ and $n_{\alpha}$ has the same order as $\alpha$. Let $N_{i,\alpha, r}'=B(\mathbb{C}_i,3r)\cap\Gamma_\alpha$ and $N_{i,\alpha, r}''=B(\mathbb{C}_i,6r)\cap\Gamma_\alpha$, we have $N_{i,\alpha, r}'\subset B(\mathbb{C}_i,3r)$ and $N_{i,\alpha, r}''\subset B(\mathbb{C}_i,6r)$, so the volumes of $N_{i,\alpha, r}'$ and $N_{i,\alpha, r}''$
are bounded by $O(r^d)$ for all $1\le i\le n_{\alpha}$. Define $S_{i,\alpha,r}=\int_{\mathbb{C}_i}V(dy)$, $S_{i,\alpha,r}'=\int_{N_{i,\alpha,r}'}V(dy)$ and $S_{i,\alpha,r}''=\int_{N_{i,\alpha,r}''}V(dy)$. Clearly, $V(dx)$ is a function of $\Xi_{B(x,r)}\cap \Gamma_\alpha$, $S_{i,\alpha,r}'$, $S_{i,\alpha,r}''$
$V_{\alpha,r}-S_{i,\alpha,r}'$ and $V_{\alpha,r}-S_{i,\alpha,r}''$ are functions of $\Xi_{B(\mathbb{C}_i,4r)}$, $\Xi_{B(\mathbb{C}_i,7r)}$, $\Xi_{\Gamma_\alpha\backslash B(\mathbb{C}_i,2r)}$ and $\Xi_{\Gamma_\alpha\backslash B(\mathbb{C}_i,5r)}$ respectively. For convenience, we write $\tilde{\Xi}$ as an independent copy of $\Xi$ and $\tilde{V}_{\alpha,r}$, $\tilde{V}(dx)$, $\tilde{S}'_{i,\alpha,r}$,  $\tilde{S}_{i,\alpha,r}''$ as the corresponding counterparts of $V_{\alpha,r}$, $V(dx)$, $S'_{i,\alpha,r}$, $S_{i,\alpha,r}''$.

For the first term in \Ref{Stein1}, 
\begin{align}
	&\mathbb{E}f'(V_{\alpha,r})\nonumber
	\\=& \mathbb{E}\(\int_{\Gamma_\alpha}V(dx)\)^2\mathbb{E}f'(V_{\alpha,r})\nonumber
	\\=& \sum_{i=1}^{n_\alpha}\mathbb{E}\(S_{i,\alpha,r}S_{i,\alpha,r}'\)\mathbb{E}f'(V_{\alpha,r})+ \mathbb{E}\(\(\int_{\Gamma_\alpha}V(dx)\)^2-\sum_{i=1}^{n_\alpha}S_{i,\alpha,r}S_{i,\alpha,r}'\)\mathbb{E}f'(V_{\alpha,r})\nonumber	
	%\\=& \sum_{i=1}^{n_\alpha}\mathbb{E}\(S_{i,\alpha,r}S_{i,\alpha,r}'\)\mathbb{E}f'(V_{\alpha,r})+ \mathbb{E}\(\(\int_{\Gamma_\alpha}V(dx)\)^2-\sum_{i=1}^{n_\alpha}S_{i,\alpha,r}S_{i,\alpha,r}'\)\mathbb{E}f'(V_{\alpha,r})\nonumber
	\\=& \sum_{i=1}^{n_\alpha}\mathbb{E}\(S_{i,\alpha,r}S_{i,\alpha,r}'\)(\mathbb{E}f'(V_{\alpha,r})-\mathbb{E}f'(V_{\alpha,r}-S_{i,\alpha,r}'')+\mathbb{E}f'(V_{\alpha,r}-S_{i,\alpha,r}''))+\epsilon_1\nonumber
	\\=& \sum_{i=1}^{n_\alpha}\(\mathbb{E}\(S_{i,\alpha,r}S_{i,\alpha,r}'\)\mathbb{E}\(\int_{0}^{S_{i,\alpha,r}''}f''(V_{\alpha,r}-x)dx\)\)+ \sum_{i=1}^{n_\alpha}\mathbb{E}\(S_{i,\alpha,r}S_{i,\alpha,r}'f'(\tilde{V}_{\alpha,r}-\tilde{S}_{i,\alpha,r}'')\)\nonumber
	\\&+\epsilon_1 \nonumber
	\\=& \sum_{i=1}^{n_\alpha}\mathbb{E}\(S_{i,\alpha,r}S_{i,\alpha,r}'\(f'(\tilde{V}_{\alpha,r}-\tilde{S}_{i,\alpha,r}'')-f'(V_{\alpha,r}-S_{i,\alpha,r}'')+f'(V_{\alpha,r}-S_{i,\alpha,r}'')\)\)+\epsilon_1+\epsilon_2\nonumber
	\\=& \sum_{i=1}^{n_\alpha}\mathbb{E}\(S_{i,\alpha,r}S_{i,\alpha,r}'f'(V_{\alpha,r}-S_{i,\alpha,r}'')\)+\epsilon_1+\epsilon_2+\epsilon_3,\label{Stein2}
\end{align} where 
\begin{align*}
	\epsilon_1=& \mathbb{E}\(\(\int_{\Gamma_\alpha}V(dx)\)^2-\sum_{i=1}^{n_\alpha}S_{i,\alpha,r}S_{i,\alpha,r}'\)\mathbb{E}f'(V_{\alpha,r});\\
	\epsilon_2=& \sum_{i=1}^{n_\alpha}\(\mathbb{E}\(S_{i,\alpha,r}S_{i,\alpha,r}'\)\mathbb{E}\(\int_{0}^{S_{i,\alpha,r}''}f''(V_{\alpha,r}-x)dx\)\) ;\\
	\epsilon_3=&  \sum_{i=1}^{n_\alpha}\mathbb{E}\(S_{i,\alpha,r}S_{i,\alpha,r}'\(f'(\tilde{V}_{\alpha,r}-\tilde{S}_{i,\alpha,r}'')-f'(V_{\alpha,r}-S_{i,\alpha,r}'')\)\).
\end{align*} 

Using the same idea as in the proof of Lemma~\ref{newlma4} and the fact that $\|f'\|\le \sqrt{\frac{2}{\pi}}$, we have 
\begin{equation}\label{e1}|\epsilon_1|\le \sqrt{\frac2{\pi}}\left|\mathbb{E}\sum_{i=1}^{n_\alpha}S_{i,\alpha,r}\(\int_{\Gamma_\alpha}V(dx)-S_{i,\alpha,r}'\)\right|\le \alpha^{-1}
\end{equation} for \qcon{} $r>C_{\thecon}\ln(\alpha)$\setcounter{cproofd}{\thecon} for some positive constant $C_{\thecon}$. 

For $\epsilon_2$, using the fact that $\|f''\|\le 2$ and Remark~\ref{moment1}, we have
\begin{align}
	|\epsilon_2|\le& 2  \sum_{i=1}^{n_\alpha}\(\mathbb{E}\(\left|S_{i,\alpha,r}S'_{i,\alpha,r}S''_{i,\alpha,r}\right|\)\)\nonumber \\\le& O(\sigma_{\alpha,r}^{-3}n_\alpha\max_{i}\{\Vol (N_{i,\alpha, r}'),\Vol (N_{i,\alpha, r}'')\}^5)=O(\alpha^{-\frac{3}{2}\nu+1}r^{5d}).\label{e2}
\end{align}

To bound $\epsilon_3$, we can use the coupling method. Since $(S_{i,\alpha,r}, S_{i,\alpha,r}')$ is a function of $\Xi_{B(\mathbb{C}_i,4r)}$ and $V_{\alpha,r}-S_{i,\alpha,r}''$ is a function of $\Xi_{\Gamma_\alpha\backslash B(\mathbb{C}_i,5r)}$, from the EDD and \cite[p.~254]{BHJ}, there is a coupling $(X_{i,1},X_{i,2},X_{i,3})$ and $(Y_{i,1},Y_{i,2},Y_{i,3})$ of $(S_{i,\alpha,r}, S_{i,\alpha,r}',f'(\tilde{V}_{\alpha,r}-\tilde{S}_{i,\alpha,r}''))$ and $(S_{i,\alpha,r}, S_{i,\alpha,r}',f'(V_{\alpha,r}-S_{i,\alpha,r}''))$ such that $\mathbb{P}(E_i):=\mathbb{P}((X_{i,1},X_{i,2},X_{i,3})\neq(Y_{i,1},Y_{i,2},Y_{i,3}))\le \alpha^{\theta_0/d} \qcon{} C_{\thecon}r^{\theta_0} \qcon{} e^{-C_{\thecon}r}$\qcon{}, for all $1\le i\le n_\alpha$. Then for $r>C_{\thecon}\ln(\alpha)$ where $C_{\thecon}$\setcounter{cproofe}{\thecon} is large enough, $\mathbb{P}(E_i)\le \qcon{} C_{\thecon}\alpha^{-6}$. Now,
\begin{align}
	|\epsilon_3|\le &  \sum_{i=1}^{n_\alpha}\left|\mathbb{E}\(S_{i,\alpha,r}S_{i,\alpha,r}'\(f'(\tilde{V}_{\alpha,r}-\tilde{S}_{i,\alpha,r}'')-f'(V_{\alpha,r}-S_{i,\alpha,r}'')\)\)\right|\nonumber
	\\=&\sum_{i=1}^{n_\alpha}\left|\mathbb{E}\(X_{i,1}X_{i,2}X_{i,3}-Y_{i,1}Y_{i,2}Y_{i,3}\)\right|\nonumber
	\\ \le &\sum_{i=1}^{n_\alpha}\(\left|\mathbb{E}\(X_{i,1}X_{i,2}X_{i,3}\mathbf{1}_{E_i}\)\right|+\left|\mathbb{E}\(Y_{i,1}Y_{i,2}Y_{i,3}\mathbf{1}_{E_i}\)\right|\). \label{xiae3}
\end{align}
However, applying H\"older's inequality, the fact that $\|f'\|\le \sqrt{\frac{2}{\pi}}\le 1$ and Remark~\ref{moment1}, we have
\begin{align*}&\left|\mathbb{E}\(X_{i,1}X_{i,2}X_{i,3}\mathbf{1}_{E_i}\)\right|\le \|X_{i,1}X_{i,2}\|_{3/2}\mathbb{P}(E_i)^{1/3}\le O\(r^{10d/3}\)\sigma_{\alpha,r}^{-2} \mathbb{P}(E_i)^{1/3},\\
&\left|\mathbb{E}\(Y_{i,1}Y_{i,2}Y_{i,3}\mathbf{1}_{E_i}\)\right|\le \|Y_{i,1}Y_{i,2}\|_{3/2}\mathbb{P}(E_i)^{1/3}\le O\(r^{10d/3}\)\sigma_{\alpha,r}^{-2} \mathbb{P}(E_i)^{1/3}.
\end{align*}
Combining with \Ref{xiae3}, we get
\begin{equation}
	|\epsilon_3|\le n_\alpha\sigma_{\alpha,r}^{-2} O(r^{10d/3})\max_{1\le i\le n_\alpha} \mathbb{P}(E_i)^{1/3}\le O(\alpha^{-1-\nu}r^{10d/3}). \label{e3}
\end{equation}

For the second term in \Ref{Stein1}, we have
\begin{align}
	&\mathbb{E}V_{\alpha,r}f(V_{\alpha,r})\nonumber
	\\=&\sum_{i=1}^{n_\alpha}\mathbb{E}\(S_{i,\alpha,r}f(V_{\alpha,r})\)\nonumber
	\\=&\sum_{i=1}^{n_\alpha}\mathbb{E}\(S_{i,\alpha,r}\(f(V_{\alpha,r}-S_{i,\alpha,r}')+f(V_{\alpha,r})-f(V_{\alpha,r}-S_{i,\alpha,r}')\)\)\nonumber
	\\=&\sum_{i=1}^{n_\alpha}\mathbb{E}\(S_{i,\alpha,r}S_{i,\alpha,r}'\int_{0}^1f'(V_{\alpha,r}-uS_{i,\alpha,r}')du\)+\sum_{i=1}^{n_\alpha}\mathbb{E}\(S_{i,\alpha,r}f(\tilde{V}_{\alpha,r}-\tilde{S}_{i,\alpha,r}')\)\nonumber
	\\&+\sum_{i=1}^{n_\alpha}\mathbb{E}\(S_{i,\alpha,r}\(f(V_{\alpha,r}-S_{i,\alpha,r}')-f(\tilde{V}_{\alpha,r}-\tilde{S}_{i,\alpha,r}')\)\)\nonumber
	\\=&\sum_{i=1}^{n_\alpha}\mathbb{E}\(S_{i,\alpha,r}S_{i,\alpha,r}'\int_{0}^1f'(V_{\alpha,r}-uS_{i,\alpha,r}')du\)+\epsilon_4+\epsilon_5,\label{Stein3}
\end{align} where 
\begin{align*}
	\epsilon_4=& \sum_{i=1}^{n_\alpha}\mathbb{E}\(S_{i,\alpha,r}f(\tilde{V}_{\alpha,r}-\tilde{S}_{i,\alpha,r}')\);\\
	\epsilon_5=& \sum_{i=1}^{n_\alpha}\mathbb{E}\(S_{i,\alpha,r}\(f(V_{\alpha,r}-S_{i,\alpha,r}')-f(\tilde{V}_{\alpha,r}-\tilde{S}_{i,\alpha,r}')\)\).
\end{align*} Since $f(\tilde{V}_{\alpha,r}-\tilde{S}_{i,\alpha,r}')$ is independent of $S_{i,\alpha,r}$ for all $1\le i\le n_\alpha$, 
\begin{align}
	\epsilon_4=\sum_{i=1}^{n_\alpha}\mathbb{E}\(S_{i,\alpha,r}\)\mathbb{E}\(f(\tilde{V}_{\alpha,r}-\tilde{S}_{i,\alpha,r}')\)=0 \label{e4}
\end{align} from the definition of $S_{i,\alpha,r}$ and $\|f\|\le 2$.

To bound $\epsilon_5$, we can use the same idea as that for bounding $\epsilon_3$, i.e., we can construct a coupling of $(S_{i,\alpha,r},f(V_{\alpha,r}-S_{i,\alpha,r}'))$ and $(S_{i,\alpha,r},f(\tilde{V}_{\alpha,r}-\tilde{S}_{i,\alpha,r}'))$. Since $S_{i,\alpha,r}$ is a function of $\Xi_{B(\mathbb{C}_i,r)}$ and $V_{\alpha,r}-S_{i,\alpha,r}'$ is a function of $\Xi_{\Gamma_\alpha\backslash B(\mathbb{C}_i,2r)}$, from the EDD and \cite[p.~254]{BHJ}, there is a coupling $((X_{i,1},X_{i,2}),(Y_{i,1},Y_{i,2}))$ of $(S_{i,\alpha,r},f(V_{\alpha,r}-S_{i,\alpha,r}'))$ and $(S_{i,\alpha,r},f(\tilde{V}_{\alpha,r}-\tilde{S}_{i,\alpha,r}'))$ such that $\mathbb{P}\((X_{i,1},X_{i,2})\neq (Y_{i,1},Y_{i,2})\)=:\mathbb{P}(E')\le \alpha^{\theta_0/d} \qcon{} C_{\thecon}r^{\theta_0} \setcounter{xiaf}{\thecon}\qcon{} e^{-C_{\thecon}r}$\setcounter{xiag}{\thecon} for all $1\le i \le n_\alpha,$ where $C_{\thexiaf},\ C_{\thexiag}\in\mathbb{R}_+$ are independent of $i$'s.  Then we follow the steps as for \Ref{e3} to get
\begin{align}
	|\epsilon_5|\le O(\alpha^{-1-\nu}r^{10d/3}), \label{e5}
\end{align}for \qcon{}$r>C_{\thecon}\ln(\alpha)$ with sufficiently large $C_{\thecon}$\setcounter{cprooff}{\thecon}.

Also, from the fact that $\|f''\|\le 2$ and Remark~\ref{moment1}, 
\begin{align}
	|\epsilon_6|:=&\left|\sum_{i=1}^{n_\alpha}\mathbb{E}\(S_{i,\alpha,r}S_{i,\alpha,r}'f'(V_{\alpha,r}-S_{i,\alpha,r}'')\)-\sum_{i=1}^{n_\alpha}\mathbb{E}\(S_{i,\alpha,r}S_{i,\alpha,r}'\int_{0}^1f'(V_{\alpha,r}-uS_{i,\alpha,r}')du\) \right|\nonumber
	\\ \le &\sum_{i=1}^{n_\alpha}\left|\mathbb{E}\(S_{i,\alpha,r}S_{i,\alpha,r}'\int_{0}^1\(f'(V_{\alpha,r}-uS_{i,\alpha,r}')-f'(V_{\alpha,r}-S_{i,\alpha,r}'')\)du\)\right|\nonumber
	\\ \le &\sum_{i=1}^{n_\alpha}\left|\mathbb{E}\(S_{i,\alpha,r}S_{i,\alpha,r}'\int_{0}^1\int_{0}^{S_{i,\alpha,r}''-uS_{i,\alpha,r}'}f''(V_{\alpha,r}-S_{i,\alpha,r}''+v)dvdu\)\right|\nonumber
	\\ \le &2\sum_{i=1}^{n_\alpha}\mathbb{E}\(\left|S_{i,\alpha,r}S_{i,\alpha,r}'\right|\(\left|S_{i,\alpha,r}''\right|+\left|S_{i,\alpha,r}'\right|\)\)\le O(\alpha^{-\frac{3}{2}\nu+1}r^{5d}).\label{e6}
\end{align}
Combining \Ref{distance1}, \Ref{Stein1}, \Ref{Stein2}, \Ref{e1}, \Ref{e2}, \Ref{e3}, \Ref{Stein3}, \Ref{e4}, \Ref{e5} and \Ref{e6}, we obtain the bound
\begin{equation}
	d_{W}\(\frac{\bar{W}_{\alpha,r}-\mu_{\alpha}}{\sigma_\alpha},\bar{Z}_{\alpha,r}\) \le O(\alpha^{-\frac{3}{2}\nu+1}r^{5d}) \label{thm2.25}
\end{equation} for \qcon{} $r>C_{\thecon}\ln(\alpha)$\setcounter{cproofg}{\thecon}, where $C_{\thecon}=\max\{C_{\thecproofd},C_{\thecproofe},C_{\thecprooff}\}$.

From \Ref{thm2a01}, taking $r=\max\{C_{\thecproofa},C_{\thecproofc},C_{\thecproofb}, C_{\thecproofg}\}\ln(\alpha)$ for large $\alpha$, together with the bounds in \Ref{thm2a02}, \Ref{thm2a05} and \Ref{thm2.25}, we have 
$$d_W\(\frac{\bar{W}_\alpha-\mu_\alpha}{\sigma_\alpha},Z\)\le \alpha^{-\frac{3}{2}\nu+1}\ln(\alpha)^{5d}.$$

(i) If $\eta$ is range-bound, then there exists an $r_1>0$ such that $\bar{W}_{\alpha,r_1}= \bar{W}_{\alpha}$ $a.s.$ for all $\alpha$, so $d_W\(\frac{\bar{W}_\alpha-\mu_\alpha}{\sigma_\alpha},Z\)= d_W\(\frac{\bar{W}_{\alpha,r_1}-\mu_{\alpha}}{\sigma_\alpha},\bar{Z}_{\alpha,r_1}\)$. Since $\eta$ is range-bound, it is also exponentially stabilising, \Ref{thm2a03} and \Ref{thm2.25} still hold, then $d_W\(\frac{\bar{W}_\alpha-\mu_\alpha}{\sigma_\alpha},Z\)=O(\alpha^{-\frac{3}{2}\nu+1}r_1^{5d})=O(\alpha^{-\frac{3}{2}\nu+1})$, completing the proof. \qed

\noindent{\it Proof of Theorem~\ref{thm2.a1}.} The statement can be shown by replacing $\bar{W}_\alpha$, $\bar{W}_{\alpha,r}$, $\bar{Z}_\alpha$, $\bar{Z}_{\alpha,r}$, $g_\alpha(x,\Xi)$ and $\bar{R}(x,\alpha)$ by their counterparts $W_\alpha$, $W_{\alpha,r}$, $Z_\alpha$, $Z_{\alpha,r}$, $g(\Xi^x)$ and $R(x)$ in the proof of Theorem~\ref{thm2.a2}. \qed

According to Theorem~\ref{thmvar}, Corollary~\ref{thm2} and Corollary~\ref{thm2a} are special cases of Theorem~\ref{thm2.a1} and Theorem~\ref{thm2.a2} respectively.
%%%%%%%%%%%%%%%%%%%%%%%%%%%%%
%%%%%%% References   %%%%%%%%
%%%%%%%%%%%%%%%%%%%%%%%%%%%%%
\newlength{\bibsep}%{0.5ex}
\def\bibfont{\small}
\normalem

\def\ac{{Academic Press}~}
\def\aap{{Adv. Appl. Prob.}~}
\def\ap{{Ann. Probab.}~}
\def\anap{{Ann. Appl. Probab.}~}
\def\eljp{{\it Electron.\ J.~Probab.\/}~} 
\def\jap{{J. Appl. Probab.}~}
\def\jrss{{J. R. Stat. Soc.}~}
\def\jtp{{J. Theor. Probab.}~}
\def\jws{{John Wiley $\&$ Sons}~}
\def\ny{{New York}~}
\def\ptrf{{Probab. Theory Related Fields}~}
\def\sp{{Springer}~}
\def\spa{{Stochastic Process. Appl.}~}
\def\spl{{Stat. Probab. Lett.}~}
\def\sv{{Springer-Verlag}~}
\def\tpa{{Theory Probab. Appl.}~}
\def\zw{{Z. Wahrsch. Verw. Gebiete}~}

\end{document}